\documentclass[12pt,a4paper]{amsart}

\usepackage{xcolor}
\usepackage{fancyhdr,mathtools}
\usepackage{appendix}
\usepackage{amssymb,amscd,amsxtra,calc}
\usepackage{mathrsfs}
\usepackage{multirow}
\usepackage{comment}
\usepackage{enumitem}
\usepackage{oplotsymbl}
\usepackage[all]{xy}
\usepackage{longtable}
\usepackage[colorlinks,linkcolor=blue,anchorcolor=red,citecolor=green,pagebackref,linktocpage]{hyperref}
\usepackage[noadjust]{cite}
\usepackage{mathabx,epsfig}

\def\acts{\ \rotatebox[origin=c]{-90}{$\circlearrowright$}\ }
\def\racts{\ \rotatebox[origin=c]{90}{$\circlearrowleft$}\ }
\def\dacts{\ \rotatebox[origin=c]{0}{$\circlearrowright$}\ }
\def\uacts{\ \rotatebox[origin=c]{180}{$\circlearrowright$}\ }

\newcommand{\plusdot}{%
  \mathbin{
    \ooalign{
      \scalebox{1.5}{$+$}
      \raisebox{0.4ex}{\kern-1.12em\scalebox{0.8}{$\bullet$}}
    }%
  }%
}

\theoremstyle{plain}
    \newtheorem{thm}{Theorem}[section]
    \newtheorem{corollary}[thm]{Corollary}

       \newtheorem{lemma}[thm]{Lemma}
           \newtheorem{theorem}[thm]{Theorem}
            
    \newtheorem{proposition}[thm]{Proposition}

\theoremstyle{definition}
 
    \newtheorem{claim}[thm]{Claim}
    \newtheorem{definition}[thm]{Definition} 
     
         \newtheorem{example}[thm]{Example}
       \newtheorem{conjecture}[thm]{Conjecture}
    \newtheorem{question}[thm]{Question}
    
    \newtheorem{notation}[thm]{Notation}
    \newtheorem*{notation*}{Notation and Terminology}
    \newtheorem{remark}[thm]{Remark}
    \newtheorem*{ack}{Acknowledgments}
\theoremstyle{remark}

\newcommand{\PP}{\mathbb{P}}

\newcommand{\Q}{\mathbb{Q}}

\newcommand{\R}{\mathbb{R}}

\newcommand{\Z}{\mathbb{Z}}

\newcommand{\Exc}{\operatorname{Exc}}

\newcommand{\id}{\operatorname{id}}

\newcommand{\NE}{\overline{\operatorname{NE}}}
\newcommand{\Nef}{\operatorname{Nef}}
\newcommand{\Nlc}{\operatorname{Nlc}}

\newcommand{\NS}{\operatorname{NS}}

\newcommand{\PE}{\operatorname{PE}}

\newcommand{\Per}{\operatorname{Per}}

\newcommand{\Sing}{\operatorname{Sing}}

\newcommand{\Supp}{\operatorname{Supp}}

\newcommand{\N}{\operatorname{N}}

\newcommand{\Sym}{\operatorname{Sym}}

\newcommand{\Pic}{\operatorname{Pic}}

\DeclareSymbolFont{extraup}{U}{zavm}{m}{n}
\DeclareMathSymbol{\varheart}{\mathalpha}{extraup}{86}
\DeclareMathSymbol{\vardiamond}{\mathalpha}{extraup}{87}

\dedicatory{Dedicated to Professor Yongnam Lee on the occasion of his 60th birthday}

\begin{document}

\title[Endomorphisms of threefolds]
{On smooth rationally connected projective threefolds of Picard number two admitting int-amplified endomorphisms}

\author{Zelong Chen, Sheng Meng and Guolei Zhong}

\address{
\textsc{School of Mathematical Sciences, East China Normal University}\endgraf
    \textsc{Shanghai 200241, China}
}
\email{zlchen@stu.ecnu.edu.cn}

\address{
\textsc{School of Mathematical Sciences, Shanghai Key Laboratory of PMMP}\endgraf 
\textsc{East China Normal University, Shanghai 200241, China}
}
\email{smeng@math.ecnu.edu.cn}

\address{
\textsc{Center for Complex Geometry, 
	Institute for Basic Science (IBS)}\endgraf
	\textsc{55 Expo-ro, Yuseong-gu, Daejeon, 34126, Republic of Korea
}}
\email{zhongguolei@u.nus.edu, guolei@ibs.re.kr}

\begin{abstract}
    We prove that a smooth rationally connected projective threefold of Picard number two is toric if and only if it admits an int-amplified endomorphism.
    As a corollary, we show that a totally invariant smooth curve of a non-isomorphic surjective endomorphism of $\mathbb{P}^3$ must be a line when it is blowup-equivariant.
\end{abstract}



\maketitle
\tableofcontents
\section{Introduction}\label{S:Intro}

We work over an algebraically closed field $\mathbf{k}$ of characteristic zero.

As one of the most historical problems in algebraic geometry, the classification of smooth projective varieties \(X\) admitting non-isomorphic surjective endomorphisms \(f\) has seen a lot of progress in recent years. 
When \(\dim X=1\), the Hurwitz formula implies that \(X\) is either a rational or an elliptic curve. 
In the case \(\dim X=2\), the structure of such (possibly singular) varieties has been completely described by Nakayama  (see~\cite{Nak20}). 

In higher dimensions, on the one hand, to exclude the non-essential case when an endomorphism is simply a product of some endomorphism with an automorphism, it is necessary to impose a constraint on \(f\). 
The appropriate condition, introduced by the second author, turns out to be 
\textit{int-amplified}.
A surjective endomorphism \(f\colon X\to X\) of a normal projective variety \(X\) is said to be int-amplified if one of the following equivalent conditions holds (see \cite[Theorem~1.1]{Men20} and \cite[Proposition~3.7]{MZ23-survey}).
\begin{itemize}
\item There exists an ample Cartier divisor \(L\) on \(X\) such that \(f^*L - L\) is ample.   
\item Every eigenvalue of \(f^*|_{\NS(X)}\) has modulus greater than one.
\item The last dynamical degree of \(f\) is strictly greater than the other dynamical degrees of \(f\); in this case,  \(f\) is also said to be \textit{\((\dim X)\)-cohomologically hyperbolic}. 
\end{itemize}

Int-amplified endomorphisms have emerged as powerful tools in the study of general endomorphisms, particularly in reducing a non-isomorphic endomorphism to a \textit{polarized} one, i.e., an endomorphism \(f\) satisfying \(f^*H\sim qH\) for some ample divisor \(H\) and for some integer \(q>1\). 
For concrete applications in various contexts, we refer the reader to \cite{Yos21, MZZ22, MZ23str, JXZ23, KT24} and the references therein.

On the other hand, the three canonical fibrations — namely, the Kodaira fibration, the Albanese map (or more specifically, the Beauville--Bogomolov--Yau decomposition), and the maximal rationally connected fibration (or more specifically, the Cao--H\"oring decomposition) — naturally facilitate the decomposition of endomorphisms via equivariant descent.  
This framework highlights rationally connected varieties as one of the central building blocks in the study of endomorphisms of higher-dimensional varieties.

In light of the aspects introduced, it becomes a fundamental problem to characterize smooth rationally connected projective varieties admitting int-amplified endomorphisms.  
In this direction, generalizing the work of Fakhruddin \cite[Question~4.4]{Fak03}, the last two authors proposed \cite[Question~1.2]{MZg23}, stated as follows.

\begin{question}[{\cite[Question 1.2]{MZg23}}]\label{que: main}
Let \(X\) be a smooth rationally connected projective variety admitting an int-amplified  endomorphism $f$. 
Is $X$ a toric variety?
\end{question}

Nakayama \cite{Nak02} confirmed Sato's conjecture that a smooth projective rational surface admitting a non-isomorphic surjective endomorphism is toric, in particular, Question \ref{que: main} has a positive answer in dimension two. 
Question \ref{que: main} also extends the following long-standing Conjecture \ref{conj: main}
of the 1980's.

\begin{conjecture}\label{conj: main}
Let $X$ be a smooth Fano variety of Picard number one. Suppose that $X$ admits a non-isomorphic surjective endomorphism $f$.
Then $X$ is a projective space.
\end{conjecture}

Conjecture \ref{conj: main} is known for threefolds \cite{ARV99,HM03}, for rational homogeneous spaces \cite{PS89}, and for hypersurfaces \cite{Bea01}. 
We refer the reader to \cite{SZ24,SZ25,KT24} and the references therein for more information and recent progress on Conjecture \ref{conj: main}.

Building on the works \cite{ARV99, HM03}, as well as the theory of the \textit{equivariant minimal model program} (EMMP) for polarized and int-amplified endomorphisms, initiated by Zhang and the second author in \cite{MZ18, Men20, MZ20}, Zhang and the last two authors established an affirmative answer to Question~\ref{que: main} for Fano threefolds \cite{MZZ22}. 

In this paper, we take a first step toward addressing Question~\ref{que: main} for arbitrary smooth rationally connected projective threefolds.  
Our main result is stated as follows.

\begin{theorem}\label{thm: main}
A smooth rationally connected projective threefold $X$ of Picard number 2 is toric if and only if it admits an int-amplified endomorphism.
\end{theorem}

\begin{proof}
    By Theorem \ref{thm: Yoshikawa}, 
    the threefold  $X$ admitting an int-amplified endomorphism is of Fano type.
    By Remark \ref{rmk: framework}, the proof of Theorem \ref{thm: main} consists of the following:
    \begin{itemize}
        \item Theorem \ref{thm: SF2}, the case when $X$ admits a Fano contraction to a surface,
        \item Theorem \ref{thm: SF1}, the case when $X$ admits a Fano contraction to a curve,
        \item Theorem \ref{thm: div-div}, the case when $X$ admits two divisorial contractions,
        \item Theorem \ref{thm: div-small}, the case when $X$ admits one divisorial contraction and one small contraction.
    \end{itemize}
    Another direction is clear; 
    see \cite[Lemma 4]{Nak02} and \cite[Proof of Theorem 1.4]{MZg23}.
\end{proof}

\begin{remark}
In fact, the last two cases in the proof of Theorem \ref{thm: main} do not occur and thus such a threefold has to be either a splitting \(\mathbb{P}^1\)-bundle over \(\mathbb{P}^2\) or a splitting \(\mathbb{P}^2\)-bundle over \(\mathbb{P}^1\).
In general, any smooth projective toric variety of Picard number two is a splitting projective bundle over a projective space, see \cite[Corollary 4.4]{Bat91}.
\end{remark}

\noindent
{\bf Main difficulties.}
Previously in \cite{MZZ22}, the analysis relied heavily on the existence of two \(K_X\)-negative extremal contractions for Fano threefolds. 
Notably, the case of small contractions did not arise, thanks to Mori's seminal work  \cite[Theorem~3.3]{Mor82}.  
In contrast, in the present paper, only one \(K_X\)-negative extremal contraction is guaranteed, and the second extremal contraction may be small.  
Even when the second contraction is log divisorial, little is known about its structure. 

Furthermore, for Fano threefolds in \cite{MZZ22}, when Fano contractions occur, we only need to deal with  $\mathbb{P}^1$-fibrations.  
In this paper, however, we must also handle with more intricate cases, including quadric surface fibrations and sextic del Pezzo fibrations.

\vskip 2mm

In what follows, we give an interesting application of Theorem~\ref{thm: main} to the following significant and challenging Conjecture~\ref{main-conj-linear}.  
The case of \(\mathbb{P}^2\) was completely solved due to Gurjar \cite{Gur03} or Forn{\ae}ss-Sibony \cite[Section 4]{FS94}. 
In the case of \(\mathbb{P}^3\), it is known that \(Y\) is linear if it is a prime divisor, thanks to the works of Nakayama–Zhang \cite[Theorem~1.5]{NZ10} and H\"oring \cite[Corollary~1.2]{Hor17}.  
For higher dimensional projective spaces  \(\mathbb{P}^n\), 
the \(Y\) is known to be linear if it is a smooth (or even with isolated singularities) hypersurface, by results of Beauville \cite[Theorem]{Bea01}, Cerveau–Lins Neto \cite{CL00}, Paranjape–Srinivas \cite[Proposition~8]{PS89}, and \cite[Corollary 1.4]{Yan23}.  
Note that all these known results address only the case where \(Y\) is a prime divisor, as the proofs primarily rely on properties of the sheaf of logarithmic differential forms.  
For higher-codimensional subvarieties, this conjecture remains intractable at present.

\begin{conjecture}\label{main-conj-linear}
Let $f$ be a surjective endomorphism of $X \coloneqq\PP^n$ with $\deg f \ge 2$.
Then any $f^{-1}$-invariant (or $f^{-1}$-periodic) closed subvariety $Y$ is a linear subspace of $X$.
\end{conjecture}

Let $Y\subseteq X=\mathbb{P}^n$ be a smooth closed subvariety.
We say $Y$ is {\it blowup-equivariant} with respect to a non-isomorphic surjective endomorphism $f\colon X\to X$  if there exists a surjective endomorphism \(\widehat{f}\colon \widehat{X}\to\widehat{X}\) such that \(f\circ\pi=\pi\circ\widehat{f}\). 
As an application of Theorem \ref{thm: main}, we show the linearity of blowup-equivariant smooth curve in $\mathbb{P}^3$.

\begin{corollary}\label{cor: main}
    Let $f\colon X\cong \mathbb{P}^3\to X$ be a surjective endomorphism with $\deg f \ge 2$.
    Let $Y\subseteq \mathbb{P}^3$ be a smooth curve such that the blowup $\pi\colon \textup{Bl}_Y X\to X$ is $f$-equivariant.
    Then $Y$ is a line.
\end{corollary}

\begin{proof}
    By Theorem \ref{thm: main}, we see that $\textup{Bl}_C X$ is toric.
    Let $T$ be the big torus acting on $\textup{Bl}_C X$.
    Then $\pi$ is also $T$-equivariant; 
    see \cite[Proposition 2.1]{Bri11}. 
    Therefore, the curve $C$ is also $T$-invariant and hence linear in $\mathbb{P}^3$.
\end{proof}

We conclude this section with the following example.
Motivated by this example, we expect in general that an $f^{-1}$-invariant curve in $\mathbb{P}^n$ with $\deg f\ge 2$ is blowup-equivariant up to a deformation of $f$.
\begin{example}
Consider the following family of endomorphisms parametrized by \(t\).
\begin{align*}
f_t\colon \mathbb{P}^3&\longrightarrow\mathbb{P}^3\\
(x_0:x_1:x_2:x_3)&\longmapsto (x_0^2+tx_1x_2:x_1^2:x_2^2:x_3^2)
\end{align*}
First, we assume that \(t\neq 0\). 
Let $L_{ij}$ be the line defined by $x_i=x_j=0$. 
Then it is straightforward to see that $L_{13}$ and $L_{01}$ are $f_t^{-1}$-invariant.
We claim that $\textup{Bl}_{L_{13}}\mathbb{P}^3\to\mathbb{P}^3$ is $f$-equivariant but $\textup{Bl}_{L_{01}}\mathbb{P}^3\to\mathbb{P}^3$ is not. 
Note that the blowup $\textup{Bl}_{L_{ij}}\mathbb{P}^3$ is the graph of the natural projection \(\mathbb{P}^3\dashrightarrow\mathbb{P}^1\) via 
$(x_0:x_1:x_2:x_3)\mapsto (x_i:x_j)$.
Hence, the blowup-equivariancy is equivalent to the descending \(f_t|_{\mathbb{P}^1}\) being well-defined.
For the former case, it is clear that \(f_t\) induces a well-defined endomorphism \(g_t\colon\mathbb{P}^1\to\mathbb{P}^1\) sending \((x_1\colon x_3)\mapsto (x_1^2\colon x_3^2)\).
For the latter case, however, it is not. 
Now, if we consider the endomorphism \(f_0\), then both blowups are \(f_0\)-equivariant with a similar argument.
\end{example}

\begin{ack}
The authors would like to thank De-Qi Zhang for the valuable discussions and suggestions. 
They also thank Andreas H{\"o}ring and Tianle Yang for the useful comments to improve the paper.  
The authors are grateful to the organizers of the ``Conference on Recent Progress in Algebraic Geometry'' held in South Korea in November 2024, where the discussions served as the foundation for this work. 
The second author is supported by the Shanghai Pilot Program for Basic Research, the Science and Technology Commission of Shanghai Municipality (No. 22DZ2229014), and the National Natural Science Foundation of China. The third author is supported by the Institute for Basic Science (IBS-R032-D1). 
\end{ack}

\section{Preliminary}\label{S:Pre}

\subsection{Varieties, divisors, and endomorphisms}

\begin{notation}
    Let $X$ be a normal projective variety.
    We use the following notation throughout this paper unless otherwise stated. 
\begin{longtable}{p{2cm} p{11cm}}

$\Pic(X)$    &the group of Cartier divisors of $X$ modulo linear equivalence~$\sim$\\
$\Pic^0(X)$    &the neutral connected component of $\Pic(X)$\\
$\NS(X)$    &$\Pic(X)/\Pic^0(X)$, the N\'eron-Severi group\\
$\N^1(X)$    &$\NS(X)\otimes_{\Z}\mathbb{R}$, the space of $\R$-Cartier divisors modulo numerical equivalence $\equiv$\\
$\N_{n-1}(X)$    &the space of Weil \(\mathbb{R}\)-divisors modulo weak numerical equivalence~$\equiv_w$\\
$\Nef(X)$   & the cone of nef classes in $\N^1(X)$\\
$\PE(X)$  & the cone of pseudo-effective classes in $\N^1(X)$\\
 $\PE_{n-1}(X)$   & the cone of pseudo-effective Weil \(\mathbb{R}\)-divisors in $\N_{n-1}(X)$\\
$\Supp D$   & the support of effective divisor $D=\sum a_i D_i$ which is $\bigcup_{a_i >0} D_i$, where the $D_i$ are prime divisors\\
$\rho(X)$    & Picard number of $X$ which is $\dim_{\R} \N^1(X)$
\end{longtable}
Here, two Weil \(\mathbb{R}\)-divisors \(P_1\) and \(P_2\) are said to be \textit{weakly numerically equivalent} (denoted by \(P_1\equiv_w P_2\)) if for any Cartier divisors \(H_1,\cdots,H_{n-1}\), one has \((P_1-P_2)\cdot H_1\cdots H_{n-1}=0\).
Note that two \(\mathbb{R}\)-Cartier divisors $D_1\equiv D_2$ if and only if $D_1\equiv_w D_2$; see \cite[Definition 2.2 and Lemma 2.3]{MZ18}.
\end{notation}

\begin{definition}\label{defn: Toric, FT,LCY}
Throughout this article, the varieties below will be heavily involved.
\begin{enumerate}
\item A normal projective variety $X$ is of \textit{Fano type}, if there is an effective Weil $\mathbb{Q}$-divisor $\Delta$ on $X$ such that the pair $(X,\Delta)$ has at worst klt singularities and $-(K_X+\Delta)$ is ample and $\mathbb{Q}$-Cartier; 
 the pair $(X,\Delta)$ is called {\it log Fano}.
If $\Delta=0$, then $X$ is the usual klt Fano variety.
We say that $X$ is {\it weak Fano} if $X$ is of Fano type and $-K_X$ is further nef (and hence big).
\item A normal projective variety $X$ is \textit{of Calabi-Yau type} (resp. \textit{of klt Calabi-Yau type}), if there is an effective Weil $\mathbb{Q}$-divisor $\Delta$ on $X$ such that the pair $(X,\Delta)$ has at worst lc (resp.~klt) singularities and $K_X+\Delta$ is $\mathbb{Q}$-Cartier and numerically trivial; the pair \((X,\Delta)\) is called \textit{log Calabi-Yau} (resp.~\textit{klt log Calabi-Yau}).
\item A normal variety $X$ of dimension $n$ is a \textit{toric variety} if $X$ contains a {\it big torus} $T=(\mathbf{k}^*)^n$ as an (affine) open dense subset such that the natural multiplication action of $T$ on itself extends to an action on the whole variety; in this case, the complement  $B\coloneqq X\backslash T$ is a divisor and the pair $(X,B)$ is said to be a \textit{toric pair} with $K_X+B\sim 0$.
\end{enumerate}
\end{definition}

\begin{definition}
    Let \(f\colon X\to X\) be a surjective endomorphism of a projective variety \(X\). 
    \begin{enumerate}
        \item 
        We say that \(f\) is \textit{$q$-polarized} if 
        $f^*H\sim qH$ for some ample Cartier divisor $H$ and integer $q>1$, or equivalently,  
        \(f^*|_{\textup{N}^1(X)}\) is diagonalizable with all the eigenvalues being of modulus \(q\) (\cite[Proposition 2.9]{MZ18}).
        \item 
        We say that  $f$ is \textit{int-amplified} if $f^*L-L$ is ample for some  ample  Cartier divisor $L$, or equivalently, all the eigenvalues of \(f^*|_{\textup{N}^1(X)}\) are of modulus greater than 1 (\cite[Theorem 1.1]{Men20}). 
        Clearly, every polarized endomorphism is int-amplified.
        \item  
        A subset \(S\subseteq X\) is called \textit{\(f^{-1}\)-invariant}
        (resp.~\textit{\(f^{-1}\)-periodic}) 
        if \(f^{-1}(S)=S\) (resp.~\(f^{-n}(S)=S\) for some $n>0$).
    \end{enumerate}
\end{definition}

By implementing the equivariant minimal model program, 
Yoshikawa has made significant advancements toward a solution to Question \ref{que: main}.
\begin{theorem}[{\cite[Corollary 1.4]{Yos21}}]\label{thm: Yoshikawa}
    A smooth rationally connected projective variety is of Fano type if it admits an int-amplified endomorphism.
\end{theorem}

\subsection{Singularities and positivities}
In this subsection, we study the intersection of boundaries for a given log canonical pair in the fibre spaces. 
Together with the canonical bundle formula and the dynamical assumption, we obtain several positivities of (relative) anti-log canonical divisors.
We begin with the following definition. 

\begin{definition}\label{def: branch}
\begin{enumerate}
\item Let \(f\colon X\to Y\) be a finite surjective morphism between normal varieties.
      The \textit{ramification divisor \(R_f\) of \(f\)} is defined by the formula 
      \(K_X=f^*K_Y+R_f\). 
      We denote by  \(B_f=f(\Supp R_f)\) the  (reduced) \textit{branch divisor} of \(f\). 
\item  Let $\pi\colon P\to Y$ be a surjective projective morphism with $Y$ being normal.
    Let $\nu\colon P^{n}\to P$ be the normalization, and $P^n\xrightarrow{\phi}\widetilde{Y}\xrightarrow{\sigma}Y$ 
    the Stein factorization of $\pi\circ\nu$.
    Then we define the {\it general branch divisor of \(\pi\)} by $B_\pi\coloneqq B_{\sigma}$.
\end{enumerate}
\end{definition}

When the base is a curve, the following lemma is an easy observation.
\begin{lemma}\label{lem: fat}
    Let $\pi\colon P \to Y$ be a surjective projective morphism of varieties with $Y$ being smooth and $\dim Y =1$. 
    Let $y\in B_\pi$.
    Suppose $P$ is smooth at generic points of $\pi^{-1}(y)$.
    Then $\pi^*y$ is not a reduced divisor.
\end{lemma}

\begin{proof}
    Suppose the contrary, that $\pi^*y$ is a reduced divisor.
    Let $\nu\colon P^n\to P$ be the normalization. 
    Let $P^n\xrightarrow{\phi}\widetilde{Y}\xrightarrow{\sigma}Y$ be the Stein factorization of $\pi\circ\nu$.
    By the assumption, $\nu$ is isomorphic at generic points of $\pi^{-1}(y)$.
    Then $\phi^*\sigma^*y=\nu^*\pi^*y$ is also a reduced divisor.
    Note that $\widetilde{Y}$ is a smooth curve.
    Then $\sigma^*y$ is reduced, contradicting the choice of \(y\).
\end{proof}

\begin{corollary}\label{cor: lc}
    Let $\pi\colon X\to Y$ be a projective morphism between smooth varieties with $\dim Y=1$.
    Let $P$ be a prime divisor on $X$ dominating $Y$ and let $y\in B_{\pi|_P}$. 
    Then $(X,P+\pi^*y)$ is not lc.
\end{corollary}

\begin{proof}
    Suppose the contrary, that $(X,P+\pi^*y)$ is lc.
    Then $\pi^*y=\pi^{-1}(y)$ is reduced.
    Since $X$ is smooth, $P\cap \pi^{-1}(y)$ is of pure codimension $2$.
    Then $P+\pi^{-1}(y)$ is a simple normal crossing divisor near any generic point $\eta$ of $P\cap \pi^{-1}(y)$.
    In particular, 
    $P$ is smooth at generic points of $(\pi|_P)^{-1}(y)$ and $(\pi|_P)^*y$ is a reduced divisor on $P$.
    However, this contradicts Lemma \ref{lem: fat}.
\end{proof}

In contrast to Corollary \ref{cor: lc}, we extend \cite[Theorem 1.6]{Men20} to pair version: Lemma \ref{lem: lc and pe} below. 
Lemmas \ref{lem: snc} and \ref{lem: non-normal invariant} are by-products of Lemma \ref{lem: lc and pe}. 
\begin{lemma}\label{lem: lc and pe}
    Let $f\colon X \to X $ be an int-amplified endomorphism on a normal projective variety $X$.
    Let $\Delta$ be a reduced effective $f^{-1}$-periodic Weil divisor.
    Then $-(K_X+\Delta)$ is a pseudo-effective Weil divisor.
    Suppose further $K_X+ \Delta$ is $\Q$-Cartier.
    Then $(X,\Delta)$ is lc and $-(K_X+\Delta)$ is \(\mathbb{Q}\)-effective.
\end{lemma}

\begin{proof}
    After iteration, we may assume $f^{-1}(\Delta)=\Delta$.
    By the ramification divisor formula, we have
    $$f^*(-(K_X+\Delta))-(-(K_X+\Delta))=R_f-(f^*\Delta-\Delta)$$
    where the right hand side is effective.
    By \cite[Theorem 3.3]{Men20}, the divisor $-(K_X+\Delta)$ is pseudo-effective.
    It is further \(\mathbb{Q}\)-effective when $K_X+ \Delta$ is $\Q$-Cartier by the same proof (after replacing $K_X$ by $K_X+\Delta$) of \cite[Theorem 6.2]{MZ22}.
    
    Suppose the contrary that \((X,\Delta)\) is not lc.
    Let $Z$ be any irreducible component of the non-lc locus $ \Nlc (X,\Delta)$ which is of dimension $d<\dim X$.
    By \cite[Theorem 1.4]{BH14} we have $\deg f|_Z=\deg f$ and then $f_* Z= (\deg f) Z$.
    Let $H$ be an ample divisor of $X$.
    Then $(f^i)^* H^d \cdot Z= (\deg f)^i  H^d \cdot Z$, which contradicts  \cite[Lemma 3.8]{Men20}.
\end{proof}

\begin{lemma}[cf.~Lemma \ref{lem: H snc}]\label{lem: snc}
    Let  $f\colon X \to X$ be an int-amplified endomorphism of a smooth projective variety \(X\).
    Suppose $A$ and $B$ are two different  $f^{-1}$-periodic prime divisors.
    Then \(A+B\) is a simple normal crossing divisor around the generic points of the intersection \(A\cap B\).
    In particular, $A|_B$ and $B|_A$ are reduced divisor.
\end{lemma}

\begin{proof}
    By Lemma \ref{lem: lc and pe}, the pair $(X,A+B)$ is lc.
    By Bertini's theorem, 
    we can take a sufficiently general smooth hyperplane section $H$ such that $(X,A+B+H)$ is lc.
    By \cite[Theorem 5.50]{KM98}, 
    the pair $(H, (A+B)|_H)$ is lc. 
    Repeating this process \(\dim X-2\) times,  
    we may assume that $X$ is a smooth projective surface. 
    For each point \(p\in A\cap B\),
    it follows from \cite[Lemma 2.29]{KM98} that both \(A\) and \(B\) are smooth along the point \(p\) and the blowup along $p$ will separate \(A\) and \(B\).
    So we have the transversality of \(A\cap B\).
\end{proof}

\begin{lemma}\label{lem: non-normal invariant}
    Let  $f\colon X \to X$ be an int-amplified endomorphism of a smooth projective variety \(X\).
    Let $D$ be a reduced $f^{-1}$-periodic divisor.
    Let $S\subseteq \Sing(D)$ be an irreducible closed subset with $\dim S=\dim X-2$.
    Then \(S\) is $f^{-1}$-periodic and $D$ is a normal crossing divisor around the generic point of $S$.
\end{lemma}

\begin{proof}
    By Lemma \ref{lem: lc and pe}, the pair $(X,D)$ is lc.
    Computing the discrepancy \cite[Lemma 2.29]{KM98}, 
    we see that the Cartier divisor $D$ has multiplicity exactly $2$ at the generic point of $S$ and $S$ is a log canonical centre of  $(X,D)$.
    It follows from \cite[Corollary 3.3]{BH14} that \(S\) is $f^{-1}$-periodic. 
    By a similar proof of Lemma \ref{lem: snc} and \cite[Theorem 4.15]{KM98}, we see that $D$ has exactly two analytic branches crossing transversally at the generic point of $S$.
\end{proof}

In the last part of this subsection, we show the following lemma, which extends \cite[Theorem 1.5]{Men20} to the relative version, and give an interesting application Corollary \ref{cor: lc trivial} to the canonical bundle formula.

\begin{lemma}\label{lem: K_X/Y pseudo-effective}
    Let $\pi\colon X \to Y$ be an equi-dimensional fibration of normal projective varieties whose fibres are all reduced and irreducible.
    Assume that we have the equivariant dynamical system
     $$\xymatrix{
        f \acts X \ar[r]^{\pi} &Y\racts g.
        }$$
    where $f,g$ are int-amplified endomorphisms.
    Let $E$ be a reduced and horizontal divisor such that $f^{-1}(E)=E$.
    Then $-K_{X/Y}-E$ is a pseudo-effective Weil divisor.
\end{lemma}

\begin{proof}
Since \(\pi\) is equi-dimensional and \(Y\) is normal, it follows that the inverse image \(\pi^{-1}(\textup{Sing}\,Y)\) of the singular locus of \(Y\) is of codimension \(\geq 2\) in \(X\).
Then for any prime divisor \(Q\) on \(Y\), we can define the pullback \(\pi^*Q\) first outside this bad locus (as the pullback of Cartier divisors) and then take the closure by the normality of \(X\). 
Besides, for any prime divisor \(P\) on \(X\), either \(P\) dominates \(Y\), or \(\pi(P)\) is a prime divisor on \(Y\). 
So for any divisor \(P\) on \(X\), we can write \(P=P^{h}+P^v\) where \(P^h\) is the summand, the support of whose irreducible component dominates \(X\) and \(P^v=P-P^h\) is the vertical part whose \(\pi\)-images are divisors on \(Y\).

By the ramification divisor formula, we have
        \begin{align*}
            K_X = f^* K_X+ R_f,~
            K_Y = g^* K_Y + R_g,~R_f=R_f^h+R_f^v.
        \end{align*}
By assumption, for any prime divisor \(Q\) on \(Y\), its pullback is also a prime divisor. 
Hence, if \(Q\) appears in \(R_g\), then \(g^*Q\) is non-reduced and thus \(f^*(P\coloneqq\pi^*Q)\) is non-reduced, too.
Therefore, \(\pi^*R_g\leq R_f^v\) by the commutative diagram and hence 
$$f^*(-(K_{X/Y}+E))-(K_{X/Y}+E)=\Delta$$
where
$\Delta\coloneqq (R_f^h-(f^*E-f^{-1}(E))+(R_f^v-\pi^*R_g)$ is effective.
By \cite[Theorem 3.3]{Men20}, 
we see that $-K_{X/Y}-E$ is a pseudo-effective Weil divisor.
\end{proof}

Lemma \ref{lem: K_X/Y pseudo-effective} fails if we drop the equi-dimensional assumption.
For example, if we take a one-point blowup \(X\to Y\cong\mathbb{P}^2\), then we can construct a polarized endomorphism \(f\) on \(X\) as in \cite[Lemma 4]{Nak02}.
In this case, \(K_{X/Y}\) is a non-zero effective divisor. 

\begin{corollary}\label{cor: lc trivial}
With the same notation as in Lemma \ref{lem: K_X/Y pseudo-effective}, assume further that \(K_X+E\) is \(\pi\)-trivial.
Then $(X,E+\pi^*P)$ is lc for any prime divisor \(P\) on \(Y\).
\end{corollary}

\begin{proof}
Since \((X,E)\) is lc and \(\pi\)-trivial, by the canonical bundle formula (see e.g. \cite[Section 3]{FG14}), we have
    $$K_X+E\sim_{\mathbb{Q}}\pi^*(K_Y+M_Y+B)$$
 where $B$ is the (effective) discriminant divisor of $(\pi,E)$ defined by \(B=\sum_{P}(1-b_P)P\) and \(M_Y\) is the \(b\)-nef (and hence pseudo-effective) \(\mathbb{Q}\)-divisor.
 Here,  \(P\) is any prime divisor and \(b_P\) is given by
 \[
 b_P=\max\left\{t\in\mathbb{Q}~|~(X,E+t\pi^*P) \textup{ is sub lc over the generic point of }P\right\}.
 \]
By Lemma \ref{lem: K_X/Y pseudo-effective}, 
\(-K_{X/Y}-E\) is pseudo-effective and hence so is \(-\pi^*(M_Y+B)\).
Therefore, \(B=0\) and \(M_Y\equiv 0\) which concludes our corollary. 
\end{proof}

\subsection{Endomorphisms on some del Pezzo surfaces}
We prepare some simple facts on the smooth quadric surface (which is isomorphic to \(\mathbb{P}^1\times\mathbb{P}^1\)) and the sextic del Pezzo surface (which is the blowup of \(\mathbb{P}^3\) along three general points).
The main task is to introduce the cross $ \plusdot $ and hexagon $ \hexago$ in Notations \ref{not: cross} and \ref{not: hexagon} for the purpose of better demonstrating the proofs of Theorems \ref{thm: dp8} and \ref{thm: dp6}.
\begin{lemma}\label{lem: f1f2}
    Let $f\colon X=\mathbb{P}^1\times \mathbb{P}^1 \to X$ be a surjective endomorphism.
    Then either $f=f_1 \times f_2$ or $f\circ \sigma = f_1 \times f_2$ 
    where $\sigma(x_1,x_2)=(x_2,x_1)$.
\end{lemma}

\begin{proof}
    Let $(a,b)\in X$.
    Denote by $H_1=\{a\}\times \mathbb{P}^1$ and $H_2=\mathbb{P}^1\times\{b\}$.
    Note that $\Nef(X)=\mathbb{R}_{\ge 0} H_1+\mathbb{R}_{\ge 0} H_2$.
    Then either $f^* H_1 \equiv q H_1$ or $\sigma^*f^*H_1\equiv q H_1$ for some $q>0$.
    For the first case, $f^*H_2\equiv q'H_2$, and the two projections of $X$ are both $f$-equivariant. So $f$ splits.
    For the second case, the same argument holds for $f\circ \sigma$. 
\end{proof}

\begin{notation}\label{not: cross}
    For a point $(a,b)\in \mathbb{P}^1\times \mathbb{P}^1$, denote by
    $$\plusdot_{(a,b)}\coloneqq \{a\}\times \mathbb{P}^1\cup\mathbb{P}^1\times\{b\}$$ 
    the {\it cross} $\plusdot\subseteq\mathbb{P}^1\times \mathbb{P}^1$ with $(a,b)$ being the anchor point of the crossing.

    For a subset $S\subseteq \mathbb{P}^1\times \mathbb{P}^1$, denote by 
    $$\plusdot_{S}\coloneqq \bigcup_{(a,b)\in S} \plusdot_{(a,b)}$$
    the union of {\it crosses} anchored in $S$.
\end{notation}

\begin{lemma}\label{lem: ample cross}
    Let $A$ be an ample Weil divisor on $\mathbb{P}^1\times \mathbb{P}^1$.
    Suppose $\Supp A\subseteq \plusdot_S$ for some finite set $S$.
    Then $\Supp A=\plusdot_{\Sing(\Supp A)}$.
\end{lemma}

\begin{proof}
    By the assumption, we have 
    $$\Supp A=(\bigcup_{i=1}^m\{a_i\}\times \mathbb{P}^1)\cup (\bigcup_{j=1}^n\mathbb{P}^1\times \{b_j\})$$
    where $m,n\ge 1$ by the ampleness of $A$.
    Then 
    $$\Supp A=\plusdot_T$$
    where $T=\bigcup_{1\le i\le m, 1\le j\le n}\{(a_i,b_j)\}$.
    Finally, note that $\Sing(\Supp A)=T$.
\end{proof}

By Lemma \ref{lem: f1f2}, the following lemma is clear.
\begin{lemma}\label{lem: inverse cross}
    Let $f\colon X=\mathbb{P}^1 \times \mathbb{P}^1 \to X$ be a surjective endomorphism.
    Let $S\subseteq X$.
    Then $f^{\pm}(\plusdot_S)=\plusdot_{f^{\pm}(S)}$.
\end{lemma}

\begin{lemma}\label{lem: branch cross}
    Let  $f\colon  X=\mathbb{P}^1 \times \mathbb{P}^1 \to X$ be an int-amplified endomorphism.
    Let $S\subseteq X$ be a finite subset such that $\sharp S=\sharp f^{-1}(S)$.
    Then $\plusdot_S\subseteq B_f=\plusdot_{\Sing(B_f)}$.
\end{lemma}

\begin{proof}
    Let $\sigma(x_1,x_2)=(x_2,x_1)$.
    Then $\sharp (f\circ\sigma)^{-1}(S) =\sharp S$ and $B_{f\circ \sigma}=B_f$.
    So we may assume $f= f_1 \times f_2$ by Lemma \ref{lem: f1f2}.
    Note that $\deg f_1$ and $\deg f_2$ are both greater than $1$ because $f$ is int-amplified.
    Then $B_f=(B_{f_1}\times \mathbb{P}^1)\cup (\mathbb{P}^1\times B_{f_2})=\plusdot_{\Sing(B_f)}$.
    
    Let $(a,b)\in S$.
    Then $\sharp f^{-1}(a,b)=1$ because $\sharp S=\sharp f^{-1}(S)$.
    So $\sharp f_1^{-1}(a)=1$ and hence $a\in B_{f_1}$.
    Similarly, $b\in B_{f_2}$.
    Therefore, $\plusdot_{(a,b)}\subseteq B_f$.
\end{proof}

\begin{notation}\label{not: hexagon}
    Let $X$ be a sextic del Pezzo surface.
    Denote by $\, \hexago_X$ the union of all the six $(-1)$-curves on $S$.
\end{notation}

\begin{lemma}\label{lem: ample hexagon}
    Let $X$ be a sextic del Pezzo surface.
    Let $A$ be a reduced ample divisor such that \(\Supp A\subseteq \hexago_X\).
    Then $A=\hexago_X$.
\end{lemma}

\begin{proof}
    Write $\hexago_X=\bigcup_{i=1}^6\ell_i$ such that $\ell_i\cap \ell_j\neq \emptyset$ if and only if $i-j\equiv\pm1 \mod 6$.
    Since $A$ is ample, $\Supp A$ is connected.
    Then we may write $A=\sum_{i=1}^m \ell_i$.
    Suppose $m<6$.
    Then $A\cdot \ell_1=0$, a contradiction.    
\end{proof}

\begin{lemma}\label{lem: hexagon}
    Let $f\colon X\to X$ be a surjective endomorphism of a sextic del Pezzo surface $X$.
    Then $f^{\pm}(\hexago_X)=\hexago_X$.
    Suppose further $f$ is non-isomorphic.
    Then $B_f=\hexago_X$.
\end{lemma}

\begin{proof}
    Note that $\hexago_X=\bigcup_{i=1}^6\ell_i$ is the union of all negative curves.
    The first argument is then obvious (cf.~\cite[Lemma 4.3]{MZ22}, \cite[Proposition 11]{Nak02}).

    Write $f^* \ell_i =q \ell_j$.
    Then $q= (\deg f)^{\frac{1}{2}}$ by the projection formula.
    So $f^* (\hexago_X)=q \, \hexago_X$.
    Since $\deg f>1$, we have  $q>1$ and hence $\hexago_X\subseteq B_f$.
    Note that $K_X+\hexago_X\sim 0$.
    By the ramification divisor formula, we have 
    $K_X+\hexago_X=f^*(K_X+\hexago_X)$ and hence $\hexago_X= B_f$.
\end{proof}

\subsection{Extremal contractions of smooth projective threefolds}
We adapt Mori's well-known results \cite{Mor82} on smooth projective threefolds to our specific case for the convenience of the reader.

\begin{theorem}\label{thm: Mor82-d}
    Let \(X\) be a smooth rationally connected projective threefold of Picard number two.
    Let \(\pi\colon X\to Y\) be an extremal \(K_X\)-negative birational contraction with the exceptional locus \(E\).
    Let \(C_\pi\) be a \(\pi\)-contracted curve whose anti-canonical degree \(-K_X\cdot C_\pi \) is minimal. 
    Then  one of the following types occur.

    \begin{itemize}[leftmargin=0mm, itemindent=5mm]
        \item \textup{[E1]} The \(\pi\) is the blowup of a smooth Fano threefold \(Y\)  along a smooth curve \(\ell\subseteq Y\) and \(E\) is a \(\mathbb{P}^1\)-bundle over \(\ell\) such that \(K_X=\pi^*K_Y+E\) and  \(K_X\cdot C_\pi=E\cdot C_\pi=-2\).
        \item \textup{[E2]} The \(\pi\) is the blowup of  a smooth Fano threefold \(Y\) along a point  such that \(K_X=\pi^*K_Y+E\) and \(E\cong\mathbb{P}^2\) with \(K_X\cdot C_\pi=2E\cdot C_\pi=-2\).
        \item \textup{[E3]} The \(\pi\) is the blowup of a factorial terminal Fano threefold \(Y\) along a nodal point \(z\) such that \(K_X=\pi^*K_Y+E\) and \(E\cong\mathbb{P}^1\times\mathbb{P}^1\) with \(K_X\cdot C_\pi=E\cdot C_\pi=-1\).
        \item \textup{[E4]} The \(\pi\) is the blowup of a factorial terminal Fano threefold \(Y\) along a cuspidal point \(z\) such that \(K_X=\pi^*K_Y+E\) and \(E\) is a singular quadric surface with \(K_X\cdot C_\pi=E\cdot C_\pi=-1\).
        \item \textup{[E5]} The \(\pi\) is the blowup of a 2-factorial terminal Fano threefold \(Y\) along a \(\frac{1}{2}(1,1,1)\)-singularity \(z\) such that \(K_X=\pi^*K_Y+\frac{1}{2}E\) and \(E\cong\mathbb{P}^2\) with \(K_X\cdot C_\pi=\frac{1}{2}E\cdot C_\pi=-1\).
    \end{itemize}

    Moreover, for the types E1 \(\sim\) E4, 
    the anti-canonical divisor \(-K_Y\) is base-point-free (\cite[Proposition 4.1]{KP24}; cf.~\cite{JR06}).
\end{theorem}

\begin{remark}
    It is worth noting that, 
    as proved in \cite[Theorems 1.7 and 1.8]{KP24}, 
    if  a smooth rationally connected projective threefold \(X\) of Picard number two admits a divisorial contraction of type E3 or E4, 
    then \(X\) itself is weak Fano.
    In our paper, 
    we try to give a self-contained proof from the dynamical viewpoints and we shall only use Theorems \ref{thm: Mor82-d} above and \ref{thm: Mor82-f} below.
\end{remark}

\begin{theorem}[{\cite[Theorem 3.5]{Mor82}}]\label{thm: Mor82-f}
    Let \(X\) be a smooth rationally connected projective  threefold of Picard number two.
    Let \(\pi\colon X\to Y\) be an extremal Fano contraction.
    Then one of the following types occur.
    \begin{itemize}[leftmargin=0mm, itemindent=5mm]
        \item \textup{[C1]}  \(\pi\) is a conic bundle over \(Y\cong\mathbb{P}^2\) with  the discriminant locus \(\Delta\neq\emptyset\).
        \item \textup{[C2]} there exists a rank two locally free sheaf \(\mathcal{E}\) on \(Y\cong\mathbb{P}^2\) such that \(\pi\) is the natural projection \(X=\mathbb{P}_Y(\mathcal{E})\to Y\) (cf.~\cite[Proposition and Corollary 1]{Ele82}).
        \item \textup{[D1]} \(\pi\) is a del Pezzo fibration over a smooth curve \(Y\) whose general fibre is a del Pezzo surface of degree \(\leq 6\).
        \item \textup{[D2]} \(\pi\) is a quadric surface fibration  over a smooth  curve \(Y\) (cf.~Definition \ref{def-quadric surface fibration}).
        \item \textup{[D3]} there exists a rank three locally free sheaf \(\mathcal{E}\) on \(Y\cong\mathbb{P}^1\) such that \(\pi\) is the projection \(X=\mathbb{P}_Y(\mathcal{E})\to Y\) (cf.~\cite[Proposition and Corollary 1]{Ele82}).
\end{itemize}
\end{theorem}

\section{Framework for setting up a proof}

Let $X$ be a smooth Fano type threefold with $\rho(X)=2$.
Say $(X,\Delta)$ is log Fano. 
With the cone theorem \cite[Theorem 3.7]{KM98} applied to $(X,\Delta)$, the two extremal rays of $\NE(X)$ are both contractible.
Therefore, we can introduce the following notation.
It will help us better elaborate on our proof in the subsequent sections.

\begin{notation}\label{not: frame}
    Denote by $\mathcal{S}$ the collection of $(X\racts f,Y\xleftarrow{\pi}X\xrightarrow{\tau} Z)$ satisfying:
    \begin{enumerate}
        \item $X$ is a smooth projective threefold of Fano type with $\rho(X)=2$,
        \item $\pi$ is the contraction of some $K_X$-negative extremal ray with $\dim Y$ minimal, 
        \item $\tau$ is the contraction of another (not necessarily $K_X$-negative) extremal ray, and
        \item $f$ is an int-amplified endomorphism of \(X\) with $\pi$ and $\tau$ being $f$-equivariant (see~\cite[Theorem 1.1]{MZ20}).
    \end{enumerate}
    Denote by $$\mathcal{S}_{\textup{Fano}}\coloneqq \{(X\racts f,Y\xleftarrow{\pi}X\xrightarrow{\tau} Z)\in \mathcal{S}\,|\, \pi \text{ is Fano}\}$$
    Denote by $$\mathcal{S}_{\textup{bir}}\coloneqq \{(X\racts f,Y\xleftarrow{\pi}X\xrightarrow{\tau} Z)\in \mathcal{S}\,|\, \pi \text{ and } \tau \text{ are birational}\}$$      
\end{notation}

\begin{lemma}\label{lem: S}
    $\mathcal{S}=\mathcal{S}_{\textup{Fano}}\cup \mathcal{S}_{\textup{bir}}$.
\end{lemma}

\begin{proof}
    Let $(X\racts f,Y\xleftarrow{\pi}X\xrightarrow{\tau} Z)\in \mathcal{S}\backslash \mathcal{S}_{\textup{Fano}}$.
    Then $\pi$ is not Fano and hence a birational morphism.
    
    Suppose the contrary that $\tau$ is not birational.
    Then $\dim X>\dim Z$.
    Let $\Delta$ be an effective $\mathbb{Q}$-Cartier Weil $\mathbb{Q}$-divisor 
    such that $(X,\Delta)$ is klt and $-(K_X+\Delta)$ is ample.
    Note that a general fibre $X_z$ is not contained in $\Supp \Delta$.
    Then there exists some curve $C\subseteq X_z$ such that $C\not\subseteq \Supp \Delta$.
    In particular, $\Delta\cdot C\ge 0$.
    Then $K_X\cdot C\le (K_X+\Delta)\cdot C<0$ and hence 
    $\tau$ is a contraction of some $K_X$-negative extremal ray with $\dim Z<\dim Y$.
    This is a contradiction to the minimality of $\dim Y$.
\end{proof}

\begin{remark}\label{rmk: framework}
    By Lemma \ref{lem: S} and Theorem \ref{thm: Yoshikawa}, it suffices to show Theorem \ref{thm: main} for elements in $\mathcal{S}_{\textup{Fano}}$ and $\mathcal{S}_{\textup{bir}}$ respectively.
\end{remark}

\begin{lemma}\label{lem: not E2}
    Let $(X\racts f,Y\xleftarrow{\pi}X\xrightarrow{\tau} Z)\in \mathcal{S}_{\textup{bir}}$.
    Then $\pi$ is not of type E2 in Theorem \ref{thm: Mor82-d}.
\end{lemma}
\begin{proof}
    Suppose the contrary that $\pi$ is of type E2.
    It follows from Theorem \ref{thm: Mor82-d} that $Y$ is a smooth Fano threefold of Picard number $1$ and $\pi$ is the blowup of $Y$ along a point.
    Note that $X$ then admits a $\mathbb{P}^1$ bundle over $\mathbb{P}^2$, a contradiction.
\end{proof}

\begin{theorem}[{\cite[Theorem 5.1]{MZZ22}}]\label{thm: MZZ5.1}
     Let $(X\racts f,Y\xleftarrow{\pi}X\xrightarrow{\tau} Z)\in \mathcal{S}$ with $X$ being Fano.
     Then $X$ is isomorphic to either $\mathbb{P}_{\mathbb{P}^2}(\mathcal{O}_{\mathbb{P}^2}\oplus \mathcal{O}_{\mathbb{P}^2}(e))$ with $e\le 0$ or a blowup of $\mathbb{P}^3$ along a line.
     In particular, we have $(X\racts f,Y\xleftarrow{\pi}X\xrightarrow{\tau} Z)\in \mathcal{S}_{\textup{Fano}}$.
\end{theorem}

\begin{proof}
    We only need to notice that the $X$ is a $\mathbb{P}^2$-bundle over $\mathbb{P}^1$ when it is the blowup of $\mathbb{P}^3$ along a line.
\end{proof}

\section{Fano contraction to a surface, 
Proof of Theorem \ref{thm: SF2}}\label{S:X->S}

In this section, we handle with $(X\racts f,Y\xleftarrow{\pi}X\xrightarrow{\tau} Z)\in \mathcal{S}_{\textup{Fano}}$ (cf.~Notation \ref{not: frame}) where $\dim Y=2$.
We aim to show:

\begin{theorem}\label{thm: SF2}
Let $(X\racts f,Y\xleftarrow{\pi}X\xrightarrow{\tau} Z)\in \mathcal{S}_{\textup{Fano}}$ with \(\dim Y=2\).
Then \(X\cong \mathbb{P}_{Y}(\mathcal{O}_{Y}\oplus\mathcal{O}_{Y}(e))\) for some $e<0$ where $Y\cong\mathbb{P}^2$.
In particular, $X$ is toric.
\end{theorem}

We first show that $\pi$ is smooth.
\begin{lemma}\label{lem: discriminant empty}
    Under the assumption of Theorem \ref{thm: SF2},  the discriminant locus $\Delta$ of $\pi$ is empty.
\end{lemma}

\begin{proof}
    It follows from Theorem \ref{thm: Mor82-f} that \(\pi\) is a conic bundle over  \(Y\cong \mathbb{P}^2\).
    Suppose the contrary, that $\Delta\neq \emptyset$.
    
    Note that $\Delta$ is a divisor.  
    The fibres $X_y$ with $y\in \Delta$ satisfy the following: 
    If $y$ is a smooth point of $\Delta$, then $X_y$ is a reducible conic in $\mathbb{P}^2$. 
    If $y$ is a double point of $\Delta$ , then $X_y$ is isomorphic
    to a double line in $\mathbb{P}^2$ (cf. \cite[Proposition 1.2]{Bea77}).
    
    By \cite[Lemmas 7.2 and 7.4]{CMZ20}, we have $g^{-1}(\Delta)=\Delta$.
    Write $\Delta=\bigcup_{i=1}^m C_i$ where $C_i$ are irreducible curves.
    By \cite[Theorem]{Gur03}, each $C_i$ is a line.
    Let $P_i\coloneqq\pi^{-1}(C_i)$.
    By the cone theorem \cite[Theorem 3.7]{KM98}, $P_i=\pi^*C_i$ is a prime divisor.
    Note that a general fibre of $\pi|_{P_i}\colon P_i\to C_i\cong \mathbb{P}^1$ has two irreducible components.
    Then $\sharp B_{\pi|_{P_i}}\ge 2$ by the Hurwitz formula (cf.~Definition \ref{def: branch}).
    In particular, $\sharp C_i\cap \Sing(\Delta)\ge 2$ for each $1\le i\le m$ and hence $m\ge 3$.
    We only need $m\ge 2$.
    Note that $P_2|_{P_1}=(\pi|_{P_1})^*(C_2|_{C_1})$ is a double line.
    So we get a contradiction to Lemma \ref{lem: snc}. 
\end{proof}

\begin{proof}[Proof of Theorem \ref{thm: SF2}]  
By Lemma \ref{lem: discriminant empty}, $\pi$ is smooth
\(\mathbb{P}^1\)-fibration.
By Theorem \ref{thm: Mor82-f}, 
there exists a rank two locally free sheaf \(\mathcal{E}\) on \(Y\) such that \(X\cong\mathbb{P}_Y(\mathcal{E})\). 
The splitting of \(\mathcal{E}\) is then guaranteed by \cite[Proposition 3]{Ame03}. 
Finally, note that $X\not\cong \mathbb{P}^2\times \mathbb{P}^1$ by Notation \ref{not: frame}.
So $e<0$.
\end{proof}

The following corollary will be served for the proof of Claim \ref{cla: div-small X+ in S_bir}.
\begin{corollary}\label{cor: fano2}
    The second contraction of $X$ in Theorem \ref{thm: SF2} is never small.
\end{corollary}
    
\begin{proof}
    Let $E\cong Y$ be a section determined by the projection $$\mathcal{O}_{Y}\oplus\mathcal{O}_{Y}(e)\to \mathcal{O}_{Y}(e).$$
    Then $\mathcal{O}_X(E)|_E\cong \mathcal{O}_{\mathbb{P}^2}(e)$.
    So $\tau$ is the $E$-negative contraction with $E$ being the exceptional divisor.  
\end{proof}

\section{Discriminant loci of quadric surface fibrations}

In this section, 
we describe the discriminant locus of a quadric surface fibration with the help of non-$\mathbb{Q}$-factorial loci.
We aim to show Theorem \ref{thm: discriminant invariant} prepared for Lemma \ref{lem: dp8}.

\subsection{Non-\texorpdfstring{$\mathbb{Q}$}{Q}-factorial loci}
We recall several facts on the 
non-$\mathbb{Q}$-factorial locus treated at the end of \cite[Section 6]{MMSZ22}.
A local ring $(R, \mathfrak{m})$ is {\it $\Q$-factorial} if for any prime ideal $\mathfrak{p}$ of height one, $\mathfrak{p} = \sqrt{(f)}$ for some $f \in p$.

\begin{definition}[{\cite[Lemma 6.20]{MMSZ22}}]
    Let $X$ be a normal variety. Denote by $\hat{\mathcal{O}}_{X,x}$ the completion of the stalk.
    Denote by 
    $$\textup{NQF}(X)=\{x\in X\,|\, \hat{\mathcal{O}}_{X,x} \text{ is not } \mathbb{Q}\text{-factorial}\}$$
    the {\it non-$\mathbb{Q}$-factorial locus} of $X$.
    Here we only consider closed points. 
    Note that $\textup{NQF}(X)\subseteq \Sing(X)$.
\end{definition}

In general, it is not clear whether $\textup{NQF}(X)$ is a closed subset of $X$ or not.
The reason we work on the completion ring is due to the following functorial result. 

\begin{lemma}[{\cite[Lemma 6.20]{MMSZ22}}]\label{lem: inverse NQF}
Let $f\colon X\to Y$ be a finite surjective morphism of normal varieties.
Then $f^{-1}(\textup{NQF}(Y))\subseteq\textup{NQF}(X) $.
\end{lemma}

From the following lemma, we know that the non-$\mathbb{Q}$-factorial locus being empty implies the usual $\mathbb{Q}$-factoriality of a variety.
\begin{lemma}[{\cite[Lemma 6.21]{MMSZ22}}]\label{lem: Q-factorial}
    Let \(X\) be a normal variety and \(x\in X\) a closed point.
    Then $\hat{\mathcal{O}}_{X,x}$ being $\mathbb{Q}$-factorial implies that $\mathcal{O}_{X,x}$ is $\mathbb{Q}$-factorial; the variety $X$ is $\mathbb{Q}$-factorial if $\mathcal{O}_{X,x}$ is $\mathbb{Q}$-factorial for every $x\in X$.
\end{lemma}

\subsection{Normal morphisms}

\begin{definition}\label{def: normal}
    Let $\pi\colon X\to Y$ be a morphism of varieties. 
    Let $x\in X$, a not necessarily closed point, and $y=\pi(x)$. 
    We say that $\pi$ is {\it normal} at $x$ if $\pi$ is flat at $x$, and the scheme-theoretical fibre $X_y$ is geometrically normal at $x$ over $k(y)$.
    We say that $\pi$ is a {\it normal morphism} if $\pi$ is normal at every point of $X$.

    Suppose $\pi$ is flat. 
    By \cite[Theorem 12.1.6]{Gro66}, 
    the set $\{y\in Y\,|\, X_y \text{ is normal}\}$ is open.
    In particular, a flat morphism $\pi$ is normal if and only if $X_y$ is normal for every closed point $y\in Y$.
\end{definition}

\begin{proposition}[{\cite[Proposition 6.8.3]{Gro65}}]\label{prop: normal}
    Normal morphisms are stable under base change and composition.
\end{proposition}

\subsection{Quadric surface fibrations}

\begin{definition}\label{def-quadric surface fibration}
    Let $\pi\colon X\to Y$ be a projective morphism between varieties.
    We say $\pi$ is a {\it quadric surface fibration} if the following conditions hold.
    \begin{enumerate}
        \item $\pi$ is flat.
        \item Each fibre of $\pi$ is isomorphic to a (possibly singular) quadric surface in $\mathbb{P}^3$.
        \item A general fibre of $\pi$ is smooth.
    \end{enumerate}
    The \textit{discriminant} \(\Delta\) of \(\pi\) is the set of points over which \(\pi\) is not smooth.
\end{definition}

Since flat morphisms are stable under base change, we also have
\begin{proposition}\label{prop: quadric base change}
    Quadric surface fibrations are stable under base change.
\end{proposition}

\begin{proposition}\label{prop: quadric normal}
    Let $\pi\colon X\to Y$ be a quadric surface fibration with $Y$ being normal.
    Then $X$ is normal.
\end{proposition}

\begin{proof}
    Note that a singular quadric surface in $\mathbb{P}^3$ is normal.
    Then $\pi$ is a normal morphism.
    By Proposition \ref{prop: normal}, we see that $X$ is normal.
\end{proof}

\begin{lemma}\label{lem: quardric smooth}
    Let $\pi\colon X\to Y$ be a quadric surface fibration.
    Suppose the restriction $\Pic_{\mathbb{Q}}(X)\to \Pic_{\mathbb{Q}}(X_y)$ is surjective for a general $y\in Y$.
    Then $\pi$ is smooth.
\end{lemma}
    
\begin{proof}
    Let $\ell$ be a line in $X_y\cong \mathbb{P}^1\times \mathbb{P}^1$ with $y$ general.
    By the assumption, 
    there exists a $\mathbb{Q}$-Cartier divisor $D$ on $X$ such that $D|_{X_y}\sim_{\mathbb{Q}} \ell$.
    Let $A$ be a $\pi$-ample divisor on $X$. 
    Then
    $$0<D\cdot A\cdot X_y=D\cdot A\cdot X_{y'},$$ 
    $$0=D^2\cdot X_y= D^2\cdot X_{y'}$$
    for any $y'\in Y$.
    
    Suppose the contrary, that $\pi$ is not smooth.
    Then $X_{y'}$ is a singular quadric surface in $\mathbb{P}^3$ for some $y'\in Y$.
    Note that $\rho(X_{y'})=1$.
    Then we have $D|_{X_{y'}}\equiv 0$ by the vanishing $D^2\cdot X_{y'}=0$.
    This contradicts $D|_{X_{y'}}\cdot A|_{X_{y'}}>0$.
\end{proof}

Now we focus on the situation when the base is a curve.
\begin{lemma}\label{lem: NQF}
    Let $\pi\colon X\to Y$ be a quadric surface fibration where $Y$ is smooth of dimension one. 
    Suppose the restriction $\Pic_{\mathbb{Q}}(X_{\textup{reg}})\to \Pic_{\mathbb{Q}}(X_y)$ is surjective for general $y$. 
    Then $\textup{NQF}(X)\cap X_y\neq\emptyset$ if $X_y$ is singular.
\end{lemma}

\begin{proof}
    By Proposition \ref{prop: quadric normal}, \(X\) is normal.    
    Suppose the contrary, that $\textup{NQF}(X)\cap X_{y'}=\emptyset$ for some singular fibre $X_{y'}$.
    Since $\dim Y=1$, there are only finitely many singular fibres.
    After shrinking $Y$, we may assume $\pi$ has exactly the singular fibre $X_{y'}$.
    Then $\textup{NQF}(X)=\emptyset$ and hence $X$ is $\mathbb{Q}$-factorial by Lemma \ref{lem: Q-factorial}.
    In particular, 
    the restriction $\Pic_{\mathbb{Q}}(X)\to \Pic_{\mathbb{Q}}(X_{\textup{reg}})$ is surjective.
    By the assumption,
    we see that $\Pic_{\mathbb{Q}}(X)\to \Pic_{\mathbb{Q}}(X_y)$ is surjective if $X_y$ smooth.
    By Lemma \ref{lem: quardric smooth}, $\pi$ is smooth, a contradiction.
\end{proof}

\begin{theorem}\label{thm: discriminant invariant}
    Let $\pi\colon X\to Y$ be a quadric surface fibration over a smooth projective curve $Y$.
    Suppose $X$ admits a surjective endomorphism $f$ such that $\pi$ is $f$-equivariant. 
    Suppose further the restriction $\Pic_{\mathbb{Q}}(X_{\textup{reg}})\to \Pic_{\mathbb{Q}}(X_y)$ is surjective for general $y$.
    Then the discriminant locus $\Delta$ of $\pi$ is $(f|_Y)^{-1}$-invariant.
\end{theorem}

\begin{proof}
    By Lemma \ref{lem: inverse NQF}, we have $f^{-1}(\textup{NQF}(X))\subseteq \textup{NQF}(X)$.
    Let $S=\overline{\textup{NQF}(X)}\subseteq \Sing(X)$.
    Note that $X$ is normal by Proposition \ref{prop: quadric normal}.
    Then $f^{-1}(S)\subseteq S$ by \cite[Lemma 7.2]{CMZ20} and hence $f^{-1}(S)=S$.
    By Lemma \ref{lem: NQF}, we have $\pi(S)=\Delta$.
    Then we are done by \cite[Lemma 7.5]{CMZ20}.
\end{proof}

\section{Fano contraction to a curve, Proof of Theorem \ref{thm: SF1}}\label{S:X->C}

In this section, we handle with $(X\racts f,Y\xleftarrow{\pi}X\xrightarrow{\tau} Z)\in \mathcal{S}_{\textup{Fano}}$ (cf.~Notation \ref{not: frame}) where $\dim Y=1$.
We aim to show:

\begin{theorem}\label{thm: SF1}
Let $(X\racts f,Y\xleftarrow{\pi}X\xrightarrow{\tau} Z)\in \mathcal{S}_{\textup{Fano}}$ with $\dim Y=1$.
Then \(\pi\) is a \(\mathbb{P}^2\)-bundle over \(Y\cong \mathbb{P}^1\); in particular, \(X\) is toric.
\end{theorem}

\begin{proof}
By \cite[Theorem 1.10]{Men20},
\(f\) descends to an int-amplified  endomorphism \(g\) on \(Y\).
By \cite[Theorem 5.1]{Fak03}, 
the set of \(g\)-periodic points $\Per(g)$ is Zariski dense in $Y$.
Let $U\subseteq Y$ be a non-empty Zariski open subset such that $X_y$ is smooth for $y\in U$.
For $y \in \Per(f)\cap U \ne \emptyset$, 
there exists $s_y>0$ such that $g^{s_y}(y)=y$ and hence $f^{s_y}(X_y)=X_y$.
Note that $f^{s_y}|_{X_y}$ is still int-amplified.
Then $X_y$ is a smooth toric del Pezzo surface by \cite[Theorem 3]{Nak02}.
By further applying Theorem \ref{thm: Mor82-f}, 
the general fibre of \(\pi\) belongs to one of the three cases: \(\mathbb{P}^2\), a smooth quadric surface $\mathbb{P}^1\times \mathbb{P}^1$, or a del Pezzo surface of degree 6.
The latter two cases are excluded in Theorems \ref{thm: dp8} and \ref{thm: dp6}.
Therefore, \(X\) is a smooth \(\mathbb{P}^2\)-bundle over \(\mathbb{P}^1\); 
in particular, by the triviality of the Brauer group of \(\mathbb{P}^1\) as well as Grothendieck's splitting theorem, 
our \(X\) is a splitting \(\mathbb{P}^2\)-bundle over \(\mathbb{P}^1\) and hence toric by \cite[(7.3.1), Page 337]{CLS11}.
\end{proof}

The following corollary will be served for the proof of Claim \ref{cla: div-small X+ in S_bir}.
\begin{corollary}\label{cor: fano1}
    Suppose further in Theorem \ref{thm: SF1} that $-K_X$ is nef and the second contraction $\tau$ is small.
    Then $\kappa(X,D)\ge 1$ for any effective divisor $D$ on $X$.
\end{corollary}

\begin{proof}
    Let $X\cong \mathbb{P}_{\mathbb{P}^1}(\mathcal{O}_{\mathbb{P}^1}\oplus\mathcal{O}_{\mathbb{P}^1}(a)\oplus\mathcal{O}_{\mathbb{P}^1}(b))$ with $0\ge a\ge b$ in Theorem \ref{thm: SF1}.
    Since $\tau$ is small, we have $b<0$.
    Let $E$ the divisor determined by the projection 
    $$\mathcal{O}_{\mathbb{P}^1}\oplus\mathcal{O}_{\mathbb{P}^1}(a)\oplus\mathcal{O}_{\mathbb{P}^1}(b)\to \mathcal{O}_{\mathbb{P}^1}(a)\oplus\mathcal{O}_{\mathbb{P}^1}(b).$$
    Let $C$ be a section determined by the projection 
    $$\mathcal{O}_{\mathbb{P}^1}\oplus\mathcal{O}_{\mathbb{P}^1}(a)\oplus\mathcal{O}_{\mathbb{P}^1}(b)\to \mathcal{O}_{\mathbb{P}^1}(b)$$
    Then $E\cdot C=b$.
    By the relative Euler sequence, we have
    $$-K_X\sim 3E-(a+b-2)F$$
    where $F$ is a fibre of $\pi$.

    Since $-K_X$ is nef, we have $-K_X\cdot C=2b-a+2\ge 0$.
    If $a=b<0$, then $E\cong \mathbb{P}^1\times \mathbb{P}^1$ and the $E$-negative curves covers $E$, a contradiction because $\tau$ is small.
    So $a>b$ and hence $(a,b)=(0,-1)$.
    Note that 
    $$h^0(X,nE)=h^0(\mathbb{P}^1,\Sym^n( \mathcal{O}_{\mathbb{P}^1}\oplus \mathcal{O}_{\mathbb{P}^1}\oplus\mathcal{O}_{\mathbb{P}^1}(-1)))=n+1.$$
    Therefore, we have $\kappa(X,E)=1$.
    Note that $\PE(X)$ is generated by $E$ and $F$ where $F$ is the fibre divisor of $\pi$ with $\kappa(X,F)=1$.
    So $\kappa(X, E)\ge 1$ for any effective divisor $E$.
\end{proof}

\subsection{Quadric surface fibration case, Proof of Theorem \ref{thm: dp8}}
\begin{theorem}\label{thm: dp8}
Under the assumption of Theorem \ref{thm: SF1}, the general fibre of \(\pi\) is not isomorphic to a smooth quadric surface.
\end{theorem}

We prepare the following lemma first.

\begin{lemma}\label{lem: dp8}
    Under the assumption of Theorem \ref{thm: SF1} with the general fibre of $\pi$ isomorphic to $\mathbb{P}^1\times \mathbb{P}^1$, we have the following.
    \begin{enumerate}
        \item The second contraction $\tau\colon X \to Z$ is birational and $X$ is not Fano.
        \item After iteration, there exists an $f^{-1}$-invariant prime divisor $P$ dominating $Y$.
        \item For any $y$ not in the discriminant locus $\Delta$ of $\pi$, the intersection set $X_y\cap P=\plusdot_{\Sing(X_y\cap P)}$  (cf. Notation \ref{not: cross}) and
        the restriction $P|_{X_y}$ is a divisor of type $(p,p)$ in $\mathbb{P}^1\times \mathbb{P}^1$ with $p=1$ or $2$.
        \item Let $\nu\colon P^n\to P$ be the normalization and $P^n\xrightarrow{\phi}\widetilde{Y}\xrightarrow{\sigma} Y$ the Stein factorization of $\pi|_P\circ \nu$.
        Then $\sharp B_{\sigma}\ge 2$ (cf.~Definition \ref{def: branch}).
        \item The discriminant locus $\Delta$ is $(f|_Y)^{-1}$-invariant.
    \end{enumerate}
\end{lemma}

\begin{proof}
    By Theorem \ref{thm: Mor82-f}, every fibre of $\pi$ is an irreducible and reduced quadric in  $\mathbb{P}^3$.
    So $\pi$ is a quadric surface fibration as we defined in Definition \ref{def-quadric surface fibration}.

    Let $y\not\in \Delta$.
    Then $X_y\cong \mathbb{P}^1\times \mathbb{P}^1$.
    By the adjunction formula, the divisor $K_X|_{X_y}=K_{X_y}$ is of type $(-2,-2)$ on $\mathbb{P}^1\times \mathbb{P}^1$.
    Let $D$ be any Cartier divisor on $X$.
    Since $\rho(X)=\rho(Y)+1$, we have $D-tK_X\in \pi^*\N^1(Y)$ for some rational number $t$.
    In particular, $D|_{X_y}$ is of type $(-2t,-2t)$ on $\mathbb{P}^1\times \mathbb{P}^1$.

    Denote by $g\coloneqq f|_Y$.
    
    \begin{enumerate}[leftmargin=*]
        \item  
        Suppose the contrary $\dim X > \dim Z$.
        Then $\tau$ is also a $K_X$-negative contraction of an extremal ray.
        By Kleiman's ampleness criterion, $-K_X$ is ample.
        By Theorem \ref{thm: MZZ5.1}, it is only possible that \(X
        \cong \mathbb{P}_{\mathbb{P}^2}(\mathcal{O}_{\mathbb{P}^2}\oplus \mathcal{O}_{\mathbb{P}^2}(e))\) for some $e\le 0$.
        However, $\mathbb{P}_{\mathbb{P}^2}(\mathcal{O}_{\mathbb{P}^2}\oplus \mathcal{O}_{\mathbb{P}^2}(e))$
        admits a section isomorphic to $\mathbb{P}^2$.
        Note that there is no fibre of $\pi$ isomorphic to $\mathbb{P}^2$. Then we obtain $\mathbb{P}^2$ dominating $Y\cong \mathbb{P}^1$, a contradiction.
        
        \item
        By (1), let $E$ be the $\tau$-exceptional locus.
        Then $f^{-1}(E)=E$ and $E$ dominates $Y$. 
        If $\dim E=2$, then $E$ is a prime divisor by \cite[Proposition 2.5]{KM98} and we simply let $P=E$.
        So we may assume $\dim E=1$ and let $C\subseteq E$ be an irreducible curve.
        After iteration, we may assume $f^{-1}(C)=C$. 
        Define 
        $$\Sigma\coloneqq\{y\in Y\,|\, \sharp (\pi|_C)^{-1}(y)<\deg \pi|_C\}.$$
        Since $Y$ is a smooth curve, we see that $\Sigma$ is a finite set.

        We claim that $g(\Delta)\subseteq \Delta$ and $g(\Sigma)\subseteq \Sigma$.
        For any $y',y\in Y$ with $g(y')=y$, we have that $f|_{X_{y'}}\colon X_{y'}\to X_y$ is surjective by \cite[Lemma 7.3]{CMZ20}.
        Note that $\rho(X_y)=2$ if $y\not\in \Delta$ and $\rho(X_y)=1$ if $y\in \Delta$.
        So $g(\Delta)\subseteq \Delta$.       
        Note that $f^{-1}(C)=C$.
        Then 
        $f|_{(\pi|_C)^{-1}(y')}\colon (\pi|_C)^{-1}(y')\to (\pi|_C)^{-1}(y)$ is also surjective.
        In particular, $\sharp (\pi|_C)^{-1}(y')\ge \sharp (\pi|_C)^{-1}(y)$ and hence $g(\Sigma)\subseteq \Sigma$.
        
        For a subset $S\subseteq \mathbb{P}^1\times \mathbb{P}^1$, recall Notation \ref{not: cross} that $\plusdot_S$ is the union of crosses with anchored crossing points in $S$.
        Construct a closed subset as follows
        \[Q \coloneqq \overline{\bigcup_{y \not \in \Delta\cup \Sigma}\plusdot_{X_y\cap C}}.\]
        Note that $Q$ contains $C$. 
        Let $y \not \in \Delta\cup\Sigma$ and  $y'\in g^{-1}(y)$.
        Then $y'\not\in \Delta\cup\Sigma$.
        By Lemma \ref{lem: inverse cross}, we have
        $$(f|_{X_{y'}})^{-1}(\plusdot_{X_y\cap C})=\plusdot_{X_{y'}\cap C}$$
        and hence
        $$f^{-1}\left(\bigcup_{y \not \in \Delta\cup \Sigma}\plusdot_{X_y\cap C}\right)\subseteq \bigcup_{y \not \in \Delta\cup \Sigma}\plusdot_{X_y\cap C}.$$
        By Lemma \cite[Lemma 7.2]{CMZ20}, we have $f^{-1}(Q)\subseteq Q$ and hence $f^{-1}(Q)=Q$.

        Let $U\subseteq Y\backslash (\Delta\cup \Sigma)$ be a non-empty open subset of $Y$ such that $X_y\not\subseteq B_f$ and $B_f|_{X_y}$ is reduced for any $y\in U$.
        Let $y\in U$ and $y'\in g^{-1}(y)$.
        Then we have $$B_{f|_{X_{y'}}}=B_f|_{X_y}.$$
        Let $A= \tau^* A_Z$ where $A_Z$ is an ample divisor on $Z$.
        As $\rho (Z)=1$, we have $(f|_Z)^*A_Z \equiv q A_Z$ for some integer $q>1$. 
        Then $f^* A \equiv qA$ and $A$ is $\pi$-ample.
        Note that $A|_{X_y}$ and $A|_{X_{y'}}$ are of the same type $(d,d)$ on $\mathbb{P}^1\times \mathbb{P}^1$ for some $d>0$.
        Moreover, we have $(f|_{X_{y'}})^*(A|_{X_y})\equiv qA|_{X_{y'}}$.
        In particular, $f|_{X_{y'}}$ induces a polarized endomorphism after we identify $X_{y'}$ with $X_y$ via an arbitrary isomorphism.
        Note that 
        $$\sharp X_{y'}\cap C=\sharp (f|_{X_{y'}})^{-1}(X_y\cap C)=\sharp X_y\cap C=\deg \pi|_C.$$
        By Lemma \ref{lem: branch cross}, we have
        $$\plusdot_{X_y\cap C}\subseteq B_{f|_{X_{y'}}}\subseteq B_f$$
        and hence
        $$Q\cap \pi^{-1}(U)\subseteq B_f.$$
        By the construction, some irreducible component \(P\) of $Q$ is a prime divisor dominating $Y$.
        After iteration, $f^{-1}(P)=P$.
        \item 
        Let $U$ be the open set in (2).
        Note that $P\cap X_y\subseteq B_{f|_{X_{y'}}}=\plusdot_{\Sing(B_{f|_{X_{y'}}})}$ for any $y\in U$ with $y'\in g^{-1}(U)$ by Lemma \ref{lem: branch cross}.
        Note that $P|_{X_y}$ is ample.
        By Lemma \ref{lem: ample cross}, 
        we have
        $X_y\cap P=\plusdot_{\Sing(X_y\cap P)}$. 
        Now let $y\in Y\backslash (\Delta\cup U)$.
        If $g^{-1}(y)=y$ after iteration, then $P|_{X_y}=X_y\cap P\subseteq B_{f|_{X_y}}$ is a reduced ample divisor (see Lemma \ref{lem: snc}), and then the same argument holds. 
        If $y$ is not $g^{-1}$-periodic, then 
        there exists some $y'\in U$ with $g^s(y')=y$ for some $s>0$.
        Note that $f^s(X_{y'}\cap P)=X_y\cap P$.
        Then $X_y\cap P=\plusdot_{\Sing(X_y\cap P)}$ by Lemma \ref{lem: inverse cross}.

        By Lemma \ref{lem: lc and pe}, $-(K_X+P)$ is pseudo-effective. 
        Let $y\not\in \Delta$.
        By the adjunction formula, 
        it follows that $-(K_{X_y}+P|_{X_y})$ is pseudo-effective.
        So $P|_{X_y}$ is a divisor of type $(p,p)$ with $p=1$ or $2$.
        \item By (3), the general fibre of $\pi|_P\colon P\to Y$ has $2p$ irreducible components.
        So $\deg \sigma =2p\ge 2$.
        Since $Y\cong \mathbb{P}^1$, we have $\sharp B_\sigma\ge 2$ by the simply connectedness of \(\mathbb{P}^1\) and $\mathbb{A}^1$.
        \item
        Let $\widetilde{X}\coloneqq X\times_Y \widetilde{Y}$ and deonte by $\widetilde{\pi}\colon\widetilde{X}\to \widetilde{Y}$ and $\widetilde{\sigma}\colon\widetilde{X}\to X$ the induced morphisms.
        Then $\widetilde{\pi}$ is again a quadric surface fibration and $\widetilde{X}$ is normal by Proposition \ref{prop: quadric normal}. 
        By the diagonal embedding $\theta \colon \widetilde{Y} \hookrightarrow \widetilde{Y}\times_{Y} \widetilde{Y}$,
        the irreducible component $\theta (\widetilde{Y})$ is a section over $\widetilde{Y}$.
        Let $\widetilde{P}_1$ be the irreducible component of $P\times_{Y} \widetilde{Y}$ dominating $\theta (\widetilde{Y})$ via the dominant rational map $P\times_{Y} \widetilde{Y}\dashrightarrow \widetilde{Y}\times_{Y} \widetilde{Y}$.
        Note that the dominant rational map $P\dashrightarrow \widetilde{Y}$ has irreducible and reduced general fibres.
        Then $\widetilde{P}_1\cap \widetilde{X}_{\widetilde{y}}$ is a line  for general $\widetilde{y}\in \widetilde{Y}$.
        In particular, the restriction map $\Pic_{\mathbb{Q}}(\widetilde{X}_{\textup{reg}})\to \Pic_{\mathbb{Q}}({\widetilde{X}_{\widetilde{y}}})$ is surjective 
        because $\Pic_{\mathbb{Q}}({\widetilde{X}_{\widetilde{y}}})=\langle \widetilde{P}_1|_{\widetilde{X}_{\widetilde{y}}}, K_{\widetilde{X}}|_{\widetilde{X}_{\widetilde{y}}}  \rangle$.
        
        Note that $f|_P$ lifts to a surjective endomorphism of $P^n$.
        By \cite[Lemma 5.2]{CMZ20}, there is a surjective endomorphism $\widetilde{g}\colon \widetilde{Y}\to \widetilde{Y}$ such that $\sigma\circ \widetilde{g}=g\circ\sigma$.
        So $f$ lifts to a surjective endomorphism $\widetilde{f}$ of $\widetilde{X}$.
        Denote by $\widetilde{\Delta}$ the discriminant locus of $\widetilde{\pi}$.
        By Theorem \ref{thm: discriminant invariant}, 
        we have $\widetilde{g}^{-1}(\widetilde{\Delta})=\widetilde{\Delta}$.
        Note that $\widetilde{\Delta}=\sigma^{-1}(\Delta)$.
        Then $g^{-1}(\Delta)=\Delta$.
    \end{enumerate}
    Lemma \ref{lem: dp8} is proved.
\end{proof}

\begin{remark}
    Lemma \ref{lem: dp8}(3) does not imply that $P|_{X_y}$ is a reduced divisor for $y\not\in \Delta$.
\end{remark}

\begin{proof}[Proof of Theorem \ref{thm: dp8}]
    We adopt Lemma \ref{lem: dp8} and use its notation there.

    We claim that $(X, P+X_{y_0})$ is lc for any $y_0\in B_{\pi|_P}\neq\emptyset$.
    Once the claim is proved, we obtain a contradiction to  Corollary \ref{cor: lc}.

\par \vskip 0.5pc 
  \noindent  
    {\bf Case: $P|_{X_y}$ is of type $(2,2)$ for $y\not\in \Delta$.}
    Then $(K_X+P)|_{X_y}\sim 0$ and hence $\pi\colon(X,P)\to Y$ is an lc-trivial fibration.
    By Corollary \ref{cor: lc trivial}, it follows that $(X, P+ X_y)$ is lc for any $y \in Y$.

\par \vskip 0.5pc 
    \noindent
    {\bf Case: $P|_{X_y}$ is of type $(1,1)$ for $y\not\in \Delta$.}
    Then $P|_{X_y}$ is a reduced divisor supported by a single cross $\plusdot_{o_y}$ with $o_y$ being the point of the crossing.
    Hence $n_P$ is isomorphic at $P_y\backslash o_y$ and $\sigma$ is unramified over $y$.
    In particular, we have $B_{\sigma}\subseteq \Delta$.
    Note that $g^{-1}(\Delta)=\Delta$ and $B_{\sigma}\neq \emptyset$.
    Let $y_0\in B_{\pi|_P}=B_{\sigma}\neq\emptyset$. 
    Then $f^{-1}(X_{y_0})=X_{y_0}$ after iteration.
    So $(X, P+X_{y_0})$ is lc by Lemma \ref{lem: lc and pe}.
\end{proof}

\subsection{Sextic del Pezzo fibration case, Proof of Theorem \ref{thm: dp6}}
\begin{theorem}\label{thm: dp6}
Under the assumption of Theorem \ref{thm: SF1}, the general fibre of \(\pi\) is not isomorphic to a sextic del Pezzo surface.
\end{theorem}

We prepare the following lemma first like Lemma \ref{lem: dp8}, but there are some differences.
We are not able to show the total invariance of the discriminant locus.
However, we are lucky to show the total invariance of the horizontal part of the branch divisor.
\begin{lemma}\label{lem: dp6}
    Under the assumption of Theorem \ref{thm: SF1} with the general fibre of $\pi$ isomorphic to a sextic del Pezzo surface, 
    we have the following.
    \begin{enumerate}
        \item Let $P$ be the $\pi$-horizontal part of $B_f$.
        Then $P$ is an $f^{-1}$-invariant prime divisor.
        \item For any $y$ not in the discriminant locus $\Delta$ of $\pi$, the divisor $P|_{X_y}=\hexago_{X_y}$ (cf. Notation \ref{not: hexagon}).
        \item Let $\nu\colon P^n\to P$ be the normalization and $P^n\xrightarrow{\phi}\widetilde{Y}\xrightarrow{\sigma} Y$ the Stein factorization of $\pi|_P\circ \nu$.
        Then $\sharp B_{\sigma}\ge 2$.
    \end{enumerate}
\end{lemma}

\begin{proof}
    Let $S$ be a sextic del Pezzo surface.
    Recall Notation \ref{not: hexagon} that $\hexago_S$ is the union of all the six $(-1)$ curves in $S$.
    Then $-K_S\sim \hexago_S$.
    Note that $X_y$ is isomorphic to $S$ if and only if $y\not\in \Delta$.
    Let $y\not\in \Delta$.
    By the adjunction formula, we have $-K_X|_{X_y}=-K_{X_y}\sim \hexago_{X_y}$.
    Let $D$ be any Cartier divisor $D$ on $X$.
    Since $\rho(X)=\rho(Y)+1$, we have $D-tK_X\in \pi^*\N^1(Y)$ for some rational number $t$.
    In particular, $D|_{X_y}\equiv -t\, \hexago_{X_y}$.

    Denote by $g\coloneqq f|_Y$.
    \begin{enumerate}[leftmargin=*]
        \item
        Let $U\subseteq Y\backslash (\Delta\cup g(\Delta))$ be a non-empty open subset of $Y$ such that $X_y\not\subseteq B_f$ and $B_f|_{X_y}$ is reduced for any $y\in U$.
        Let $y\in U$ and $y'\in g^{-1}(y)$.
        The restriction $f|_{X_{y'}}\colon X_{y'}\to X_y$ is surjective by \cite[Lemma 7.3]{CMZ20}.
        Since $f$ is int-amplified, there exists a $\pi$-ample divisor $A$ such that $f^*A\equiv qA$ for some integer $q>1$.
        Moreover, we have $(f|_{X_{y'}})^*(A|_{X_y})\equiv qA|_{X_{y'}}$.
        In particular, $f|_{X_{y'}}$ induces a polarized endomorphism after we identify $X_{y'}$ with $X_y$ via an arbitrary isomorphism.
        By Lemma \ref{lem: hexagon}, we have 
        $$(f|_{X_{y'}})^{-1}(\hexago_{X_y})=\hexago_{X_{y'}}$$ and $$\hexago_{X_y}=B_{f|_{X_{y'}}}=P\cap X_y$$
        where $P$ is the $\pi$-horizontal part of $B_f$.
        Then 
        $$f^{-1}(P\cap \pi^{-1}(U))\cap \pi^{-1}(U)\subseteq P$$
        and hence $f^{-1}(P)=P$, see \cite[Lemma 7.2]{CMZ20}.

        Let $P_1$ be an irreducible component of $P$.
        Note that $P_1$ dominates $Y$.
        Then $P_1$ is $\pi$-ample.
        Let $y\in U$.
        Then $P_1|_{X_y}\subseteq \hexago_{X_y}$ is a reduced ample divisor.
        By Lemma \ref{lem: ample hexagon}, we have $P_1|_{X_y}=\hexago_{X_y}$.
        So $P=P_1$ is a prime divisor.
        \item 
        Let $y\not\in \Delta$.
        Note that $P|_{X_y}\sim \hexago_{X_y}$.
        So it suffices to show $X_y\cap P=\hexago_{X_y}$.
        If $g^{-1}(y)=y$ after iteration, then $X_y\cap P\subseteq B_{f|_{X_y}}=\hexago_{X_y}$ by Lemma \ref{lem: hexagon}.
        If $y$ is not $g^{-1}$-periodic, there exists some $y'\in U$ with $g^s(y')=y$ for some $s>0$.
        Note that $f^s(X_{y'}\cap P)=X_y\cap P$.
        Then $X_y\cap P=\hexago_{X_y}$ by Lemma \ref{lem: hexagon}.
        \item Note that $X_y\cap P$ has six irreducible components.
        So $\deg \sigma=6$ and hence $\sharp B_{\sigma}\ge 2$ by the simply  connectedness of \(\mathbb{P}^1\) and $\mathbb{A}^1$.
    \end{enumerate}
    Lemma \ref{lem: dp6} is proved.
\end{proof}

\begin{proof}[Proof of Theorem \ref{thm: dp6}]
    We adopt Lemma \ref{lem: dp6} and use its notation there.

    Note that $K_X+P$ is $\pi$-trivial.
    So $K_X+P \in \pi^* \Pic(Y)$,
    and hence $\pi\colon (X,P)\to Y$ is an lc-trivial fibration by Lemma \ref{lem: lc and pe}.
    By Corollary \ref{cor: lc trivial}, it follows that $(X, P+ X_y)$ is lc for any $y \in Y$.
    Let $y_0\in B_{\pi|_P}=B_{\sigma}\neq\emptyset$.
    Then we obtain a contradiction to Corollary  \ref{cor: lc}.
\end{proof}

\section{Two divisorial contractions, Proof of Theorem \ref{thm: div-div}}

In this section, we handle with $(X\racts f,Y\xleftarrow{\pi}X\xrightarrow{\tau} Z)\in \mathcal{S}_{\textup{bir}}$ (cf.~Notation \ref{not: frame}) where $\tau$ is a divisorial contraction.
We aim to show:

\begin{theorem}\label{thm: div-div}
There is no $(X\racts f,Y\xleftarrow{\pi}X\xrightarrow{\tau} Z)\in \mathcal{S}_{\textup{bir}}$ with $\tau$ being divisorial.
\end{theorem}

\begin{proof}
    Suppose the contrary that there exists $(X\racts f,Y\xleftarrow{\pi}X\xrightarrow{\tau} Z)\in \mathcal{S}_{\textup{bir}}$ with $\tau$ being divisorial.
    Denote by $E$ (resp. $F$) the $\pi$-exceptional (resp.~$\tau$-exceptional) divisor.
    Note that the $\tau$-exceptional curve must be $E$-positive; otherwise, the divisor $-E$ would be nef, which is absurd.
    In particular, we have $E\neq F$ by \cite[Lemma 3.39]{KM98} and
    the pseudo-effective cone $\PE(X)$ is generated by $E$ and $F$.
    Note that $-K_X$ is big because $X$ is of Fano type.
    So there exist coprime positive integers $x_1,x_2,x_3$ such that $$x_1K_X+x_2E+x_3F\equiv 0.$$
    We try to derive a contradiction based on this.

    \begin{claim}\label{cla: EF1}
       The intersection  $E\cap F$ dominates both $\pi(E)$ and $\tau(F)$.
    \end{claim}

    \begin{proof}
        If $E\cap F$ does not dominate $\tau(F)$, 
        then $E$ is $\tau$-trivial and hence 
        $E= \tau^*\tau(E)$ by the cone theorem \cite[Theorem 3.7]{KM98}.
        Note that $\tau(E)$ is an ample \(\mathbb{Q}\)-Cartier divisor because $\rho(Z)=1$ and \(Z\) is \(\mathbb{Q}\)-factorial (see \cite[Proposition 3.36]{KM98}). 
        So $E$ is nef, a contradiction to  \cite[Lemma 3.39]{KM98} again.
        Similarly, we have that $E\cap F$ dominates $\pi(E)$.
    \end{proof}

    \begin{claim}\label{cla: EF2}
        We have $\dim \pi(E)=0$ and $\dim \tau(F)=1$.
        Moreover, the contraction $\pi$ is of type E3, E4 or E5.
    \end{claim}

    \begin{proof}
        By Lemma \ref{lem: not E2}, the contraction $\pi$ is not of type E2.
        
        Suppose $\pi$ is of type E1.
        It follows from Theorem \ref{thm: Mor82-d} that $Y$ is a smooth Fano threefold of Picard number $1$ and $\pi$ is the blowup of $Y$ along a smooth curve $\pi(E)$.
        Then \(Y\cong \mathbb{P}^3\) by \cite{ARV99} or \cite{HM03}.
        Denote by $g\coloneqq f|_Y$.
        As \(\pi(F)\) is a \(g^{-1}\)-invariant prime divisor, 
        by \cite{Hor17} and \cite[Section 5]{NZ10}, it is linear. 
        Note that $E\cap F$ dominates $\pi(E)$ by Claim \ref{cla: EF1}.
        It follows that \(\pi(F)\) contains a $g^{-1}$-invariant curve \(\pi(E)\).
        By \cite[Theorem]{Gur03}, 
        we see that \(\pi(E)\) is linear as well.
        Then \(X\) is the blowup of $\mathbb{P}^3$ along a line and the second contraction $\tau$ is a $\mathbb{P}^2$-bundle over $\mathbb{P}^1$, a contradiction.

        Suppose $\tau(F)$ is a point.
        Then \(\pi(E)\) has to be a curve and thus \(\pi\) is of type E1, 
        for otherwise, 
        \(\pi\) and \(\tau\) would contract some common curve. 
        However, this is impossible as we argued above.
    \end{proof}

    By Claim \ref{cla: EF2} and Theorem \ref{thm: Mor82-d}, 
    the divisor $E$ is isomorphic to either $\mathbb{P}^2$ or an irreducible (possibly singular) quadric surface.
    Let $C_{\pi}$ be a line in $E$.
    Then $K_X\cdot C_{\pi}=-1$.

    Recall that $f$ is $q$-polarized (cf.~\cite[Corollary 3.12]{MZ18}).
    Denote by $\Delta_E$ (resp.~$\Delta_F$, $\Delta_{E+F}$) the divisor $R_f-(q-1)E$ (resp.~$R_f-(q-1)F$, $R_f-(q-1)(E+F)$).
    \begin{claim}\label{cla: div-div big}
        The divisors $-(K_X+E+F)$ and $\Delta_{E+F}$ are big.
    \end{claim}

    \begin{proof}
     By the logarithmic ramification divisor formula,  
        we have
        $$f^*(-(K_X+E+F))-(-(K_X+E+F))=\Delta_{E+F}.$$
    Hence, $\Delta_{E+F}$ is an effective divisor with $E\not\subseteq \Supp \Delta_{E+F}$ and $F\not\subseteq \Supp \Delta_{E+F}$. 
        By \cite[Theorem 1.1]{Men20}, we see that $-(K_X+E+F)$ is big if $\Delta_{E+F}$ is big.
        
        Suppose the contrary that $\Delta_{E+F}$ is not big.
        Then either $\Delta_{E+F}\sim_{\Q} \lambda E$ or $\Delta_{E+F}\sim_{\Q} \lambda F$ for some integer $\lambda\ge 0$.
        Note that $\kappa(X,E)=\kappa(X,F)=0$.
        This forces $\lambda = 0$ and thus 
        $$K_X+E+F=f^*(K_X+E+F).$$
        By \cite[Theorem 1.1]{Men20}, we have $K_X+E+F\equiv 0$.
        From \cite[Theorem 1.5]{MZg23}, we then obtain a toric pair $(X,E+F)$.
        However, the boundary divisor of a smooth projective toric threefold should have at least $4$ irreducible components (cf.~\cite[Remark 4.6]{MZ19}).
        So we get a contradiction.
    \end{proof}

    \begin{claim}\label{cla: div-div x1<x2 x1<x3}
        We have that $x_1<x_2$ and $x_1<x_3$.
    \end{claim}
    \begin{proof}
        Suppose $x_1\ge x_2$.
        Then we have the reconfigured equation
        $$(x_1-x_2)K_X+x_2(K_X+E)+x_3F\equiv 0$$
        and its $f$-pullback
        $$(x_1-x_2)(K_X-R_f)+x_2(K_X+E-\Delta_E)+x_3qF\equiv 0.$$
        Subtracting the above two equations from each other, 
        we have
        $$(x_1-x_2)R_f+x_2\Delta_E\equiv x_3(q-1)F.$$
        This contradicts Claim \ref{cla: div-div big}, noting that \(x_2>0\), \(\Delta_E>\Delta_{E+F}\) is big, and \(F\) is an extremal ray in \(\PE(X)\). 
        So we have $x_1<x_2$ and similarly $x_1<x_3$.
    \end{proof}
    
    \begin{claim}\label{cla: div-div KC>=0}
        We have $K_X\cdot C_{\tau}\ge0$ for any curve $C_{\tau}$ contracted by $\tau$.
    \end{claim}
    \begin{proof}
        Suppose the contrary that $K_X\cdot C_{\tau}<0$.
        Then $-K_X$ is ample and we get a contradiction to  Theorem \ref{thm: MZZ5.1}.
    \end{proof}

    \begin{claim}\label{cla: div-div EC>1}
        Let $C_{\tau}$ be a curve contracted by $\tau$.
        Then $E\cdot C_{\tau}>K_X\cdot C_{\tau}$ and $E\cdot C_{\tau}>1$.
    \end{claim} 

    \begin{proof}
        By Claims \ref{cla: EF2}, \ref{cla: div-div KC>=0} and Theorem \ref{thm: Mor82-d}, we have
        $$0\le K_X\cdot C_{\tau}\le K_Y\cdot\pi_*C_{\tau}+E\cdot C_{\tau}<E\cdot C_{\tau}.$$
        
        Suppose the contrary that $E\cdot C_{\tau}=1$.
        Then $K_X\cdot C_{\tau}=0$.
        The equation
        $$0=(x_1K_X+x_2E+x_3F)\cdot C_{\tau}=x_2+x_3F\cdot C_{\tau}$$ implies
        $\frac{x_2}{x_3}\in\mathbb{Z}$ 
        and hence the equation 
        $$0=(\frac{x_1}{x_3}K_X+\frac{x_2}{x_3}E+F)\cdot C_{\pi}=-\frac{x_1}{x_3}+(\frac{x_2}{x_3}E+F)\cdot C_{\pi}$$ 
        implies
        $\frac{x_1}{x_3}\in\mathbb{Z}$.
        So $x_1\ge x_3$, a contradiction to Claim \ref{cla: div-div x1<x2 x1<x3}.
    \end{proof}

    \begin{claim}\label{cla: div-div EC=2}
        There exists a curve $C_{\tau}$ contracted by $\tau$ such that $K_X\cdot C_{\tau}=0$, $E\cdot C_{\tau}=2$ and $F\cdot C_{\tau}=-2$.
    \end{claim}
    \begin{proof}
        Denote by $h\coloneqq f|_Z$.
        By Claims \ref{cla: EF1} and \ref{cla: EF2}, 
        we see that \(\tau(F)\) is a (possibly singular) curve and $E\cap F$ dominates $\tau(F)$.
        Denote by $S$ the union of $1$-dimensional components of $\Sing(F)$.
        By Lemma \ref{lem: non-normal invariant}, we have $f^{-1}(S)=S$ after iteration.
        By Lemma \ref{lem: snc}, the divisor $E+F$ has simple normal crossing near the generic points of $E\cap F$.
        So there exists a finite subset $T\subseteq F$ such that $E|_{F\backslash T}$ is a smooth divisor.
        By the generic smoothness,
        there exists a non-empty open subset 
        $U\subseteq \tau(F)\backslash\tau(T)$ such that
        \begin{itemize}
            \item $\Sing(F_z)=S\cap F_z$,
            \item $E\cap S\cap F_z=\emptyset$, and
            \item $E|_{F_z}$ is smooth, i.e., the sum of reduced points
        \end{itemize}
        for any $z\in U$.
        By \cite[Theorem 5.1]{Fak03}, 
        the set of periodic points 
        $\Per(h|_{\tau(F)})$ is Zariski dense in $\tau(F)$.
        Take $z\in \Per(h|_{\tau(F)})\cap U$ and $C_{\tau}$ any irreducible component of $F_z$.
        After iteration, we may assume $h(z)=z$ and $f(C_{\tau})=C_{\tau}$.
        Note that $C_{\tau}$ admits at most two $(f|_{C_{\tau}})^{-1}$-periodic points.
        Then $$\sharp E\cap C_\tau+\sharp S\cap C_\tau\le 2.$$
        Note that $E\cdot C_{\tau}=\sharp E\cap C_{\tau}\ge 2$ by Claim \ref{cla: div-div EC>1}.
        It follows that $E\cdot C_{\tau}=2$ and $S\cap C_\tau=\emptyset$.
        In particular, $F_z=X_z$ is a smooth rational curve and $F$ is smooth around $F_z$.
        
        Suppose $K_X\cdot C_{\tau}\ge 1$.
        Then $K_X\cdot C_{\tau}=1$ by Claim \ref{cla: div-div EC>1}.
        Recall that $K_X\cdot C_{\pi}=-1$.
        If $E\cdot C_{\pi}=-2$, then $E\equiv 2K_X$, 
        a contradiction to the bigness of \(-K_X\). 
        By Claim \ref{cla: EF2} and Theorem \ref{thm: Mor82-d}, 
        we have $E\cdot C_{\pi}=-1$.      
        Intersecting $x_1K_X+x_2E+x_3F\equiv 0$ with $C_\pi$ and $C_{\tau}$, we have equations:
        \[\left\{
            \begin{array}{rl}
                x_1+x_2&=(F\cdot C_{\pi})x_3\\
                x_1+2x_2&=-(F\cdot C_{\tau})x_3
            \end{array}
        \right.
        \]
        In particular, we obtain 
        $$x_1=(2F\cdot C_{\pi}+F\cdot C_{\tau})x_3$$
        and hence $x_3\mid x_1$, 
        a contradiction to Claim \ref{cla: div-div x1<x2 x1<x3}.
        So $K_X\cdot C_{\tau}=0$.
        
        By the adjunction formula, 
        we have 
        $$(K_X+F)\cdot C_{\tau}=K_F\cdot C_{\tau}=-2$$
        where these intersections make sense because we calculate them inside the smooth locus of $F$.
        So $F\cdot C_{\tau}= -2$.
    \end{proof}

    \noindent
    {\bf End of the proof of Theorem \ref{thm: div-div}.}
    By Claim \ref{cla: div-div EC=2}, we have equations:
    \[\left\{
            \begin{array}{rl}
                x_1-(E\cdot C_{\pi})x_2&=(F\cdot C_{\pi})x_3\\
                2x_2&=2x_3
            \end{array}
        \right.
    \]
    Then $x_1=(F\cdot C_{\pi}+E\cdot C_{\pi})x_2$.
    This contradicts Claim \ref{cla: div-div x1<x2 x1<x3}.
\end{proof}

\section{One divisorial and one small contraction}

In this section, we handle with $(X\racts f,Y\xleftarrow{\pi}X\xrightarrow{\tau} Z)\in \mathcal{S}_{\textup{bir}}$ (cf.~Notation \ref{not: frame}) where $\tau$ is a small contraction.
We aim to show:

\begin{theorem}\label{thm: div-small}
There is no $(X\racts f,Y\xleftarrow{\pi}X\xrightarrow{\tau} Z)\in \mathcal{S}_{\textup{bir}}$ with $\tau$ being a small contraction.
\end{theorem}

To prove Theorem \ref{thm: div-small}, 
we shall stick to the following notation till the end of this section.

\begin{notation}\label{not: diagram}
    Let $(X\racts f=f_X,Y\xleftarrow{\pi}X\xrightarrow{\tau} Z)\in \mathcal{S}_{\textup{bir}}$ with $\tau$ being small.
    Then we have the following diagram of equivariant dynamical systems:
    \begin{displaymath}
        \xymatrix{
            & \save[]+<-1.2pc,1.2pc>*{f_W} \restore \save[]+<0pc,1.1pc>*{\uacts} \restore W \ar[dl]_{p^+} \ar[dr]^{p} \ar@/^1pc/[drrr]^\sigma & & \\
            \save[]+<-2.4pc,0pc>*{f_{X^+}} \restore   \save[]+<-1.2pc,0pc>*{\acts}\restore X^+ \ar[rd]_{\tau^+}  &  & \save[]+<0.9pc,1.1pc>*{\rotatebox[origin=c]{150}{$\circlearrowright$} f_X} \restore X \ar@{-->}[ll]_\phi \ar[rr]_\pi \ar[ld]^{\tau} & & Y \racts f_Y \\
            &  \save[]+<-1pc,-1pc>*{f_Z} \restore \save[]+<0pc,-1.1pc>*{\dacts} \restore Z& & &}
    \end{displaymath}
    where 
    \begin{itemize}
        \item the map $\phi\colon X\dashrightarrow X^+$ is the log flip induced by the log flipping contraction $\tau$,
        \item the existence of $f_{X^+}$ follows from \cite[Lemma 6.6]{MZ18}, 
        \item the variety $W$ is the normalization of the graph of $\phi$ with two projections $p\colon W \to X$ and $p^+\colon W \to X^+$, and
        \item all the surjective endomorphism are $q$-polarized by \cite[Corollary 3.12]{MZ18}, noting that $\rho(Y)=1$.
    \end{itemize} 
    Moreover, we introduce the following notation as well as fundamental facts.
    \begin{enumerate}
        \item Denote by $E$ the $\pi$-exceptional divisor and by $C_\tau$ a fixed $\tau$-exceptional irreducible curve.
        \item Let $E_W \coloneqq p^{-1}(E)=p^*(E)$ which is Cartier and denote by $F_W$ a prime divisor contained in $p^{-1}(C_\tau)$ which dominates \(C_\tau\). 
        After iteration, we may assume that both \(E_W\) and \(F_W\) are 
         $f_W^{-1}$-invariant.
        \item Denote by $T_{f_Y}$ (resp.~\(T_{f_W}\)) the union of $f_Y^{-1}$-periodic (resp.~\(f_W^{-1}\)-periodic) prime divisors. 
        By \cite[Corollary 3.8]{MZ20}, both \(T_{f_W}\) and \(T_{f_Y}\) have only finitely many prime divisors.  
        After iteration, we may assume every $f_Y^{-1}$-periodic (resp.~\(f_W^{-1}\)-periodic) prime divisor is $f_Y^{-1}$-invariant (resp.~\(f_W^{-1}\)-invariant).
        In this way, we may assume 
        $$(f_-)^*|_{\N^1(-)}=q\id_{\N^1(-)}$$ 
        where $-$ takes $X,X^+,Y,Z,W$.
    \end{enumerate}
\end{notation}

\subsection{Intersection of exceptional loci}

In this subsection, we prove the following theorem.

\begin{theorem}\label{thm: div-small EC<=2}
    Under Notation \ref{not: diagram},
    we have the intersection $E\cdot C_{\tau}\le 2$ and $\sharp E\cap C_{\tau}=E\cdot C_{\tau}$ if $C_{\tau}$ is smooth.
\end{theorem}

\begin{lemma}\label{lem: klt lcy1}
Under Notation \ref{not: diagram}, for each \(n\in\mathbb{Z}_+\), there exists some  Weil $\mathbb{Q}$-divisor $D_{Y,n}$ such that the pair  $(Y,(1-\frac{1}{n})T_{f_Y}+D_{Y,n})$ is klt log Calabi-Yau (cf.~Definition \ref{defn: Toric, FT,LCY}).
\end{lemma}

\begin{proof}
    Recall that Theorem \ref{thm: Mor82-d} implies \(Y\) is \(\mathbb{Q}\)-factorial.
    By Lemma \ref{lem: lc and pe}, 
    the pair $(Y,T_{f_Y})$ is lc and \(-(K_Y+T_{f_Y})\) is pseudo-effective.
    Hence, 
    $(Y,(1-\frac{1}{n})T_{f_Y})$ is klt and \(-(K_Y+(1-\frac{1}{n})T_{f_Y})\) is ample, 
    noting that \(\rho(Y)=1\).
    Let $A\equiv -m(K_Y+(1-\frac{1}{n})T_{f_Y})$ be an integral very ample divisor for a sufficiently divisible \(m=m(n)\).
    It follows from Bertini's theorem that  
    a general member $H \in |A|$ intersects $T_{f_Y}$ transversally. 
    Therefore, 
    the  pair $(Y,(1-\frac{1}{n})T_{f_Y}+\frac{1}{m}H)$ is klt log Calabi-Yau.
\end{proof}

\begin{notation}\label{not: pair}
    We follow Notation \ref{not: diagram}.
    By Lemma \ref{lem: klt lcy1}, we have a series of klt log Calabi-Yau pairs $(Y, (1-\frac{1}{n})T_{f_Y}+D_{Y,n})$ for each $n \in \Z^+$.
     Let $D_{W,n}$ be the Weil $\mathbb{Q}$-divisor on $W$ obtained by 
    $$
    K_W + (1-\frac{1}{n})T_{f_W} + D_{W,n} = \sigma^*(K_Y + (1-\frac{1}{n})T_{f_Y} + D_{Y,n})\equiv 0.
    $$
    The following construction is a special application of $f$-pairs as in Yoshikawa \cite[Proposition 6.2]{Yos21}.
    For each $i\ge 0$, denote by
    \begin{align*}
        D_{Y,n,i} \coloneqq & \dfrac{1}{(\deg f_Y) ^i} (f_Y^i)_* (R_{f_Y^i} +(1-\frac{1}{n}) T_{f_Y} + D_{Y,n}) - T_{f_Y}, \\
        D_{W,n,i} \coloneqq & \dfrac{1}{(\deg f_W)^i} (f_W^i)_* (R_{f_W^i} + (1-\frac{1}{n})T_{f_W} + D_{W,n}) - (1-\frac{1}{n})T_{f_W}.
    \end{align*}
\end{notation}

\begin{lemma}\label{lem: klt lcy2}
    Under Notation \ref{not: pair}, 
    For any \(n\in\mathbb{N}\), 
    there exists \(i_n\gg 1\) such that the pair \((W, (1-\frac{1}{n})T_{f_W}+D_{W,n,i_n})\) is klt log Calabi-Yau.
\end{lemma}

\begin{proof}
    By \cite[Lemma 1.1]{FG12} and its proof, the pair $(Y, (1-\frac{1}{n})T_{f_Y}+ D_{Y,n,i})$  is klt log Calabi-Yau.
    \begin{claim}\label{cla: triviality}
        $K_W+(1-\frac{1}{n})T_{f_W}+D_{W,n,i}=\sigma^*(K_Y+(1-\frac{1}{n})T_{f_Y}+D_{Y,n,i}) \equiv 0$.
    \end{claim}
    \begin{proof}
        Recall that $f_W^*|_{\N^1(W)}=q\id_{\N^1(W)}$.
        Then 
        \begin{align*}
            (f_W^i)^* (K_W+(1-\frac{1}{n})T_{f_W}+D_{W,n,i})& 
            \equiv q^iK_W+R_{f_W^i}+(1-\frac{1}{n})T_{f_W}+D_{W,n} \\
            &=(q^i-1)K_W+R_{f_W^i}\equiv0.
        \end{align*}
        This implies that $K_W + (1-\frac{1}{n})T_{f_W}+ D_{W,n,i} \equiv 0$. 
        Therefore, we have 
        $$\Delta \coloneqq K_W + (1-\frac{1}{n})T_{f_W}+D_{W,n,i}-\sigma^*(K_Y+(1-\frac{1}{n})T_{f_Y}+D_{Y,n,i}) \equiv 0.$$
        Note that $\Supp \Delta \subseteq \Exc(\sigma)$ since $\sigma_* \Delta =0$. 
        So $ \Delta=0$ by \cite[Lemma 3.39]{KM98}.
    \end{proof}
    \begin{claim}\label{cla: effectivity}
      There exists some \(i_n\ge 0\) such that \(D_{W,n,i_n}\) is effective.
    \end{claim}
    \begin{proof}
        We refer to \cite[Claim 6.3]{Yos21} for the original proof.
        Since our pair here is slightly different, we provide a proof for the reader's convenience.

        Note that $D_{Y,n,i}$ is effective and $\tau_* D_{W,n,i} = D_{Y,n,i}$ for each $i>0$.
        Hence, it suffices to show that $\textup{mult}_{P} D_{W, i} \geq 0$ for $i \gg 1$, where $P$ is an arbitrary exceptional prime divisor of $\sigma\colon W\to Y$.
        We may assume $P$ is $f_W^{-1}$-invariant after iterstion. 
        Write $f_W^*P = rP$ and $(f_W)_*P= eP$.
        Since $(Y, (1-\frac{1}{n})T_{f_Y} + D_{Y,n})$ is klt, we have 
        \(a \coloneqq \textup{mult}_P ((1-\frac{1}{n})T_{f_W} + D_{W,n}) < 1\) 
        and 
        \begin{align*}
            \textup{mult}_{P} D_{W,n,i} = & \textup{mult}_{P} \left(\frac{{f_W^i}_*(R_{f_W^i}+(1-\frac{1}{n})T_{f_W}+D_{W,n})}{\deg f_W^i}-(1-\frac{1}{n})T_{f_W}\right)\\
            = & \frac{1}{r^i}(r^i-1+a)-(1-\frac{1}{n})\textup{mult}_{P} T_{f_W}\\
            = & \frac{1}{n}+\frac{a-1}{r^i} \ge 0
        \end{align*}
        for $i\gg1$.
        Our claim is thus proved.
    \end{proof}  
    Back to the proof of Lemma \ref{lem: klt lcy2}, 
    from Claims \ref{cla: effectivity} and \ref{cla: triviality}, 
    we obtain a klt log Calabi-Yau pair $(W, (1-\frac{1}{n})T_{f_W}+D_{W,n,i_n})$, 
    which concludes our proof.
\end{proof}

\begin{lemma}\label{lem: H}
    Under Notation \ref{not: diagram},
    there exists a hypersurface $H$ of $W$ such that the following hold. 
    \begin{enumerate}
        \item The hypersurface $H$ does not contain any irreducible component of $E_W\cap F_W$ and $\Sing(W)$.
        \item The pair $(H, E_H + F_H)$ is lc where $E_H=E_W\cap H=E_W|_H$ is Cartier and $F_H=F_W \cap H$.
    \end{enumerate}
\end{lemma}

\begin{proof}
    Note that there are only finitely many irreducible components of $E_W \cap F_W$ and $\Sing(W)$.
    So we can choose a general hypersurface $H$ such that 
    $H$ does not pass through any generic points of those irreducible components which gives (1).
    
    By \cite[Theorem 7]{Sei50}, 
    we can further choose $H$ general such that $H$ is normal.   
    Let  $U_W= W \setminus (\Sing (W) \cup (E_W \cap F_W))$ and $U_H=H \cap U_W$.
    Note that both $U_W$ and $U_H$ are smooth. 

    We follow Notation \ref{not: pair} and apply Lemma \ref{lem: klt lcy2}.
    Let $i_n\gg 1$.
    Since $E_W, F_W \subseteq T_{f_Y}$, we write
    $$(1-\frac{1}{n})T_{f_W} + D_{W,n,i_n}=(1-\frac{1}{n})(E_W + F_W) + \Delta_{W,n}$$ 
    where $\Delta_{W,n}$ is an effective Weil $\mathbb{Q}$-divisor .
    Note that $(W, (1-\frac{1}{n})(E_W+F_W)+\Delta_{W,n})$ is klt log Calabi-Yau. 
    Then $(W, (1-\frac{1}{n})(E_W+F_W) + \Delta_{W,n} + H)$ is lc by \cite[Lemma 5.17(2)]{KM98}.
    By \cite[Theorem]{Kaw07}, for each \(n\), we have
        $$K_{H}+B_n=
        (K_{W}+(1-\frac{1}{n})(E_W+F_W)+\Delta_{W,n}+H)|_{H}, $$
    where $B_n$ is the Shokurov's different of $(1-\frac{1}{n})(E_W+F_W)+\Delta_{W,n}$ such that $(H,B_n)$ is lc. 
    By the (smooth) adjunction formula for $U_H\subseteq U_W$, we have
    $$K_{U_H}+(1-\frac{1}{n})(E_{U_H}+F_{U_H})+\Delta_{U_H,n}=
    (K_{U_W}+(1-\frac{1}{n})(E_{U_W}+F_{U_W})+\Delta_{U_W,n}+H_{U_W})|_{U_H} 
    $$
    where $\Delta_{U_-,n}\coloneqq \Delta_{W,n}|_{U_-}$, $E_{U_-}\coloneqq E_W|_{U_-}$ and similarly for $F_{U_-}, H_{U_W}$.
    
    The above two formulas imply that $B_n|_{U_H}=(1-\frac{1}{n})(E_{U_H}+F_{U_H})+\Delta_{U_H,n}$.
    Note that $\dim W\backslash U_W \le 1$ and $\dim H \backslash U_H \le 0$ by (1). 
    Then $B_n$ is the extension divisor of $(1-\frac{1}{n})(E_{U_H}+F_{U_H})+\Delta_{U_H,n}$.
    Since $E_{U_H}$ and $F_{U_H}$ are integral divisors, they are reduced divisors by noting that $B_n \le 1$.
    Therefore, $E_H=E_W\cap H=E_W|_H$ and $F_H=F_W\cap H$ are the extension divisors of $E_{U_H}$ and $F_{U_H}$, respectively.
    In particular, $(1-\frac{1}{n})(E_H+F_H)\le B_n$.
    So $(H,(1-\frac{1}{n})(E_H+F_H))$ is numerically lc and hence lc; see \cite[Proposition 3.5 (2)]{Fuj12}.
    By taking $n \to \infty$, we conclude that $(H, E_H+F_H)$ is lc.
\end{proof}

\begin{remark}
    We are not able to apply Lemma \ref{lem: lc and pe} directly because  $K_W+E_W+F_W$ is not known to be $\mathbb{Q}$-Cartier.
\end{remark}

The following lemma is a slight extension of Lemma \ref{lem: snc}.
\begin{lemma}\label{lem: H snc}
    Under the notation in Lemma \ref{lem: H}, the pair
    $(H, E_H+F_H)$ is simple normal crossing near $E_H\cap F_H$.
    In particular, the restriction divisor $E_H|_{F_H}$ is reduced.
\end{lemma}

\begin{proof}
    Let $e \in E_H \cap F_H$ be a point.
    Note that $(H, E_H)$ is numerically dlt (and hence dlt) at $E_H \cap F_H$ (cf. \cite[Corollary 4.2 and Proposition 4.11]{KM98}).
    By Lemma \ref{lem: H}, $E_H$ is Cartier.
    By \cite[Corollary 5.56]{KM98}, $H$ is smooth near $E_H \cap F_H$.
    By \cite[Theorem 4.15(1)]{KM98}, $E_H+F_H$ has two exactly two analytic branches at $e$.
    In particular, $E_H$ and $F_H$ are smooth at $e$.
    Suppose $E_H$ intersects $F_H$ not transversally at $e$.
    By blowing up $H$ along $e$, we obtain an lc pair with the boundary divisor a union of three curves passing through one point, a contradiction.
\end{proof}

\begin{proof}[Proof of Theorem \ref{thm: div-small EC<=2}]
    We follow Lemma \ref{lem: H} and its notation therein. 
    Consider the following commutative diagram:
    \begin{displaymath}
        \xymatrix{
            \widetilde{F_H}\ar[r]^{n_{F_H}} \ar[d]_{\widetilde{p_F}} & F_H \ar[d]^{p_{_F}} \\
            \save[]+<-2pc,0pc>*{f_{\widetilde{C_\tau}} \acts} \restore \widetilde{C_\tau} \ar[r]^{n_{C_\tau}} & C_\tau \save[]+<2pc,0pc>*{\racts f_{C_\tau}} \restore}           
    \end{displaymath}
    where \begin{itemize}
        \item $n_{F_H}$ and $n_{C_\tau}$ are normalizations of $F_H$ and $C_\tau$ respectively.
        \item $p_F$ is the restriction of $p$ on $F_H$ and induces $\widetilde{p_F}$.
        \item $f_{C_\tau}$ is the restriction of $f_{X}\colon X \to X$ on $C_\tau$ and induces $f_{\widetilde{C_\tau}}$.
    \end{itemize}
    Since $H|_{F_W}$ is ample in $F_W$ with support $F_H=F_W \cap H$, we have
    $p_{_F}$ is surjective.
    Observe that 
    $\Supp n_{C_\tau}^*(E |_{C_\tau})=n_{C_\tau}^{-1}(E \cap C_\tau)$
    is $(f_{\widetilde{C_\tau}})^{-1}$-invariant and there are at most two $(f_{\widetilde{C_\tau}})^{-1}$-periodic points on $\widetilde{C_\tau}$, so $\sharp \Supp n_{C_\tau}^*( E|_{C_\tau}) \le 2$.

    Suppose $E \cdot C_\tau \ge 3$.
    Then there exists a point $e_1 \in \widetilde{C_\tau}$ such that 
    $\textup{mult}_{e_1} n_{C_\tau}^*( E|_{C_\tau}) \ge 2.$
    Let $e_0=n_{C_\tau}(e_1) \in E \cap C_\tau$ and choose $e_3 \in \widetilde{F_H}$ such that $\widetilde{p_F}(e_3)=e_1$.
    Let $e_2 = n_{F_H} (e_3) \in E_H \cap F_H$.
    By Lemma \ref{lem: H snc}, $H$ is smooth at $e_2$ and hence $n_{F_H} $ is isomorphic at $e_3$.
    So we have 
    $$\textup{mult}_{e_2} E_H|_{F_H}=\textup{mult}_{e_3}n_{F_H}^* (E_H|_{F_H}) = \textup{mult}_{e_3} (n_{C_\tau}\circ {\widetilde{p_F}})^*(E|_{C_\tau}) \ge 2,$$
    a contradiction to Lemma \ref{lem: H snc}.

    Suppose $C_{\tau}$ is smooth and $\sharp E\cap C_{\tau}<E\cdot C_{\tau}$.
    Then $n_{C_{\tau}}$ is an  identity map and there exists a point $e_1 \in \widetilde{C_\tau}=C_{\tau}$ such that 
    $\textup{mult}_{e_1}  E|_{C_\tau} \ge 2$. 
    Then we obtain the same contradiction.
\end{proof}

\subsection{Proof of Theorem \ref{thm: div-small}}
We first prepare the following general lemma.
\begin{lemma}\label{lem: DC>=2}
Let \(X\) be a projective variety and \(D\) a base-point-free Cartier divisor on \(X\).
For any singular curve \(C\) on \(X\), we then have \(D\cdot C\ge 2\).
\end{lemma}
\begin{proof}
Consider the morphism \(\iota\colon X\to\mathbb{P}^N\) induced by the linear system \(|D|\) such that \(D=\iota^*H\) where \(H\) is a hyperplane in \(\mathbb{P}^N\). 
Let \(W\) be the image \(\iota(X)\). 
If \(\deg \iota|_C\geq 2\), then by the projection formula, we have \(D\cdot C=(\deg \iota|_C)H|_W\cdot \iota(C)\geq 2\).
If \(\deg \iota|_C=1\), 
then \(C\) is birational onto \(\iota(C)\) and hence \(\iota(C)\) is singular as well, 
in particular, 
\(\iota(C)\) is non-linear in \(\mathbb{P}^N\) and thus 
\(D\cdot C=H|_W\cdot \iota(C)=H\cdot\iota(C)\geq 2\).
\end{proof}

In the remaining subsection, we prove Theorem \ref{thm: div-small}. 
\begin{proof}[Proof of Theorem \ref{thm: div-small}]
    Our first task is to find another totally invariant prime divisor $F$.
    \begin{claim}\label{cla:E345}
        We have $K_X\cdot C_\tau\geq 0$ and $\dim \pi(E)=0$.
        Moreover, the contraction $\pi$ is of type E3, E4 or E5.
    \end{claim}

    \begin{proof}
        By the same proof of Claim \ref{cla: div-div KC>=0}, 
        we have $K_X\cdot C_{\tau}\ge 0$.
        By Lemma \ref{lem: not E2}, \(\pi\) is not of type \(E2\).

        Suppose $\pi$ is of type E1.
        Recall the proof of Claim \ref{cla: EF2}, we see that $Y\cong \mathbb{P}^3$.
        By Theorems \ref{thm: Mor82-d} and \ref{thm: div-small EC<=2}, 
        we have
        $$0\le K_X\cdot C_{\tau}= K_Y\cdot\pi_*C_{\tau}+E\cdot C_{\tau}\le -2.$$
        This is a contradiction.
    \end{proof}

    \begin{claim}\label{cla: div-small KC=0}
        We have that $K_X\cdot C_{\tau}=0$ and hence $X$ is weak Fano.
    \end{claim} 

    \begin{proof}
        By Claim \ref{cla:E345}, \(K_X\cdot C_\tau\ge0\) and \(\pi\) is of type E3, E4 or E5. 

        Suppose $\pi$ is of type E5.
        Then we have
        $$0\leq K_X\cdot C_{\tau}= K_Y\cdot\pi_*C_{\tau}+\frac{1}{2}E\cdot C_{\tau}<1.$$
        Therefore, \(K_X\cdot C_\tau=0\).
    
        Suppose $\pi$ is of type E3 or E4.
        Then we have
        $$K_X\cdot C_{\tau}= K_Y\cdot\pi_*C_{\tau}+E\cdot C_{\tau}.$$
        If $E\cdot C_{\tau}=1$, then we still have $K_X\cdot C_{\tau}=0$.
        So we may assume $E\cdot C_{\tau}=2$.
        By Theorem \ref{thm: div-small EC<=2}, either $\sharp E\cap C_{\tau}=2$ or $C_{\tau}$ is singular.
        Note that $\pi(E)$ is a point.
        It follows that $\pi(C_{\tau})$ is a singular curve.
        By Theorem \ref{thm: Mor82-d}, the divisor $-K_Y$ is base point free.
        By Lemma \ref{lem: DC>=2}, we have $K_Y\cdot\pi_*C_{\tau}\le -2$.
        So $K_X\cdot C_{\tau}\le 0$.        
    \end{proof}

    By Claim \ref{cla: div-small KC=0}, 
    we see that $\tau$ is a $K_X$-trivial small contraction and hence it follows from \cite[Theorem 1.4.15]{IP99} (cf.~\cite[Theorem 2.4]{Kol89}, \cite{Rei83}) that there exists a flop \(X\dashrightarrow X^+\) over \(Z\) such that $X^+$ is smooth, and 
    $K_{X^+}$ is $\tau^+$-trivial and not pseudo-effective. 
    Then we have the following diagram:
    \[
    \xymatrix{
    X\ar@{-->}[rr]^{\phi}\ar[d]_{\pi}\ar[dr]^{\tau}&&X^+\ar[d]^{\pi^+}\ar[dl]_{\tau^+}\\
    Y&Z&Y^+}
    \]
    where $\pi^+\colon X^+\to Y^+$ is the $K_{X^+}$-negative contraction. 
    Note that the big divisor $-K_{X^+}$ is $\tau^+$-trivial and $\pi^+$-ample.
    So $X^+$ is again weak Fano.

    \begin{claim}\label{cla: div-small X+ in S_bir}
        We have $(X^+\racts f^+,Y^+\xleftarrow{\pi^+}X^+\xrightarrow{\tau^+} Z)\in \mathcal{S}_{\textup{bir}}$ with $\tau^+$ being small.
    \end{claim}

    \begin{proof}
        Note that $\tau^+$ is small.  
        Suppose the contrary that $\pi^+$ is a Fano contraction.
        By Corollary \ref{cor: fano2}, we have $\dim Y^+=1$.
        Note that $\kappa(X^+,\phi(E))=\kappa(X,E)=0$.
        This is a contradiction to Corollary \ref{cor: fano1}.
    \end{proof}

    By Claim \ref{cla: div-small X+ in S_bir} and Theorem \ref{thm: Mor82-d}, 
    we see that $\pi^+$ is still a divisorial contraction.
    Now we introduce the following four totally invariant divisors.
    \vskip1mm
    \begin{tabular}{@{}ll}
        $\bullet$ $E$   & the $\pi$-exceptional divisor \\
        $\bullet$ $E^+$ & the $\phi$-image $\phi(E)$ \\
        $\bullet$ $F^+$ & the $\pi^+$-exceptional divisor \\
        $\bullet$ $F$   & the $\phi$-preimage $\phi^{-1}(F^+)$ 
    \end{tabular}
    
    Note that $E^+$ is $\tau^+$-negative and $F^+$ is $\pi^+$-negative.
    Then $\PE(X)=\phi^*(\PE(X^+))$ is generated by $E$ and $F$.
    So there exists coprime positive integers $x_1, x_2, x_3$ such that
    \[\left\{
            \begin{array}{rl}
                x_1K_X+x_2E+x_3F&\equiv 0\\
                x_1K_{X^+}+x_3F^++x_2E^+&\equiv 0
            \end{array}
        \right.
    \]
    We try to derive a contradiction based on this.
    The strategy is a bit different with that in the proof of Theorem \ref{thm: div-div}.
    The trick here is that the restriction on $x_1,x_2,x_3$ is coming from both $X$ and $X^+$ symmetrically.

    Note that $f$ is $q$-polarized.
    Denote by $\Delta_E$ (resp.~$\Delta_F$, $\Delta_{E+F}$) the divisor $R_f-(q-1)E$ (resp.~$R_f-(q-1)F$, $R_f-(q-1)(E+F)$).
    By the same proofs of Claims \ref{cla: div-div big}, \ref{cla: div-div x1<x2 x1<x3} and \ref{cla: div-div EC>1}, we have the following three claims. 
    Note that Claim \ref{cla: div-small big} is based on $\kappa(X,E)=\kappa(X,F)=0$,
    Claim \ref{cla: div-small x1<x2 x1<x3} is based on Claim \ref{cla: div-small big}, and Claim \ref{cla: div-small EC=2} is based on Claim \ref{cla: div-small x1<x2 x1<x3}.
    \begin{claim}\label{cla: div-small big}
        The divisors $-(K_X+E+F)$ and $\Delta_{E+F}$ are big.
    \end{claim}

    \begin{claim}\label{cla: div-small x1<x2 x1<x3}
        We have that $x_1<x_2$ and $x_1<x_3$.
    \end{claim}

    \begin{claim}\label{cla: div-small EC=2}
        We have that $E\cdot C_{\tau}\ge 2$.
    \end{claim} 

    The above induces further restriction.
    \begin{claim}\label{cla: div-small x1=1,x2>2,x3=2}
        We have that $x_1=1$, $x_2>2$ and $x_3=2$.
    \end{claim}
    
    \begin{proof}
        By Claim \ref{cla: div-small x1<x2 x1<x3},
        we have $x_3\ge x_1+1\ge 2$ and $x_2\ge 2$.
        By Theorem \ref{thm: div-small EC<=2} and Claim \ref{cla: div-small EC=2}, 
        the equation 
        $$0=(x_1K_X+x_2E+x_3F)\cdot C_{\tau}=2x_2+(F\cdot C_{\tau})x_3$$
        implies that $\frac{2x_2}{x_3}\in \mathbb{Z}$ 
        and hence the equation
        $$0=(\frac{2x_1}{x_3}K_X+\frac{2x_2}{x_3}E+2F)\cdot C_{\pi}=-\frac{2x_1}{x_3}+(\frac{2x_2}{x_3}E+2F)\cdot C_{\pi}$$
        implies that $\frac{2x_1}{x_3}\in \mathbb{Z}$. 
        By Claim \ref{cla: div-small x1<x2 x1<x3} again, we have $\frac{2x_1}{x_3}<2$ which forces $2x_1=x_3$, and hence $\frac{x_2}{x_1}=\frac{2x_2}{x_3}\in\mathbb{Z}$.
        By the coprime assumption, we have $x_1=1$ and thus $x_3=2$ and $x_2>1$.
        Suppose $x_2=2$.
        Then we have
        $$0=(K_X+2E+2F)\cdot C_{\pi}=-1+2(E+F)\cdot C_{\pi}.$$
        This is a contradiction.
    \end{proof}

    \noindent
    {\bf End of the proof of Theorem \ref{thm: div-small}.}
    By Claim \ref{cla: div-small x1=1,x2>2,x3=2}, 
    we have 
    $$K_X+x_2E+2F\equiv 0.$$
    Replacing $(X\racts f,Y\xleftarrow{\pi}X\xrightarrow{\tau} Z)$ by $(X^+\racts f^+,Y^+\xleftarrow{\pi^+}X\xrightarrow{\tau^+} Z)$,
    we also have 
    $$K_{X^+}+x_2F^++2E^+\equiv 0$$ and hence 
    $$K_X+2E+x_2F=\phi^*(K_{X^+}+x_2F^++2E^+)\equiv 0.$$
    Therefore $x_2=2$, contradicting  Claim \ref{cla: div-small x1=1,x2>2,x3=2}.
\end{proof}

\bibliographystyle{amsalpha}

\bibliography{bib-ref}

\end{document}